\documentclass[12pt]{amsart}

\usepackage{setspace}
\usepackage{amssymb}   
\usepackage{amsmath}
\usepackage{mathrsfs}
\usepackage{stmaryrd}
\usepackage{amsthm}
\usepackage{newlfont}
\usepackage{amscd}
\usepackage{mathtools}
\usepackage{bm}
\usepackage{tikz-cd}
\usepackage{graphicx}
\usepackage{hyperref}
\usepackage{enumitem}
\usepackage[all]{xy}
\usepackage{textcomp}
\usepackage{amsbsy}

\newtheorem{thm}{Theorem}[section]
\newtheorem{thm-defn}[thm]{Theorem/Definition}
\newtheorem{lem}[thm]{Lemma}
\newtheorem{prop}[thm]{Proposition}
\newtheorem{cor}[thm]{Corollary}

\theoremstyle{definition}
\newtheorem{defn}[thm]{Definition}

\newtheorem{eg}[thm]{Example}

\theoremstyle{remark}
\newtheorem{rem}[thm]{Remark}

\numberwithin{equation}{section}

\begin{document}

\pagenumbering{arabic}

\title[Strongly divisible lattices in the imperfect residue field case]{Strongly divisible lattices and crystalline cohomology in the imperfect residue field case}
\author{Yong Suk Moon}
\date{}

\maketitle

\begin{abstract}
Let $k$ be a perfect field of characteristic $p \geq 3$, and let $K$ be a finite totally ramified extension of $K_0 = W(k)[p^{-1}]$. Let $L_0$ be a complete discrete valuation field over $K_0$ whose residue field has a finite $p$-basis, and let $L = L_0\otimes_{K_0} K$. For $0 \leq r \leq p-2$, we classify $\mathbf{Z}_p$-lattices of semistable representations of $\mathrm{Gal}(\overline{L}/L)$ with Hodge-Tate weights in $[0, r]$ by strongly divisible lattices. This generalizes the result of \cite{liu-semistable-lattice-breuil}. Moreover, if $\mathcal{X}$ is a proper smooth formal scheme over $\mathcal{O}_L$, we give a cohomological description of the strongly divisible lattice associated to $H^i_{\text{\'et}}(\mathcal{X}_{\overline{L}}, \mathbf{Z}_p)$ for $i \leq p-2$, under the assumption that the crystalline cohomology of the special fiber of $\mathcal{X}$ is torsion-free in degrees $i$ and $i+1$. This generalizes a result in \cite{cais-liu-breuil-kisin-mod-cryst-cohom}. 	
\end{abstract}

\tableofcontents

\section{Introduction} \label{sec:introduction}

Let $k$ be a perfect field of characteristic $p \geq 3$. Let $K$ be a finite totally ramified extension of $K_0 = W(k)[p^{-1}]$, and denote by $\mathcal{O}_K$ its ring of integers. In this article, we consider a complete discrete valuation field $L$ over $K$ whose reside field has a finite $p$-basis, and study lattices of semistable representations of $G_L\coloneqq \mathrm{Gal}(\overline{L}/L)$. More precisely, let $k'$ be a field extension of $k$ having a finite $p$-basis, and let $\mathcal{O}_{L_0}$ be a Cohen ring for $k'$. Let $\mathcal{O}_L = \mathcal{O}_{L_0}\otimes_{W(k)} \mathcal{O}_K$, and denote $L_0 = \mathcal{O}_{L_0}[p^{-1}]$ and $L = \mathcal{O}_L[p^{-1}]$. For semistable representations of $G_L$, it is proved in \cite{gao-integral-padic-hodge-imperfect} (cf. also \cite{brinon-crys-rep-imperfect-residue}) that weakly admissibility implies admissibility (which generalizes the result for the perfect residue field case in \cite{colmez-fontaine}), thereby giving a complete classification of semistable $\mathbf{Q}_p$-representations of $G_L$.

It is natural to ask further for a classification of $\mathbf{Z}_p$-lattices of semistable representations of $G_L$, since such an integral (or torsion) theory is very useful for studying semistable deformation ring as in the work of Liu on Fontaine's conjecture \cite{liu-fontaineconjecture}. Let $r < p-1$ be a non-negative integer. In the case of perfect residue field, Breuil conjectured such a classification of lattices of semistable representations with Hodge-Tate weights in $[0, r]$ by strongly divisible lattices \cite{breuil-integral-padic-hodge}, which is proved by Liu in \cite{liu-semistable-lattice-breuil}.

We generalize the result in \cite{liu-semistable-lattice-breuil} and classify lattices of semistable representations of $G_L$ with Hodge-Tate weights in $[0, r]$. Let $E(u) \in W(k)[u]$ be the Eisenstein polynomial for a uniformizer $\pi \in \mathcal{O}_K$, and let $S$ be the $p$-adically complete divided power envelope of $\mathcal{O}_{L_0}[\![u]\!]$ with respect to $(E(u))$. A \textit{strongly divisible lattice} $\mathscr{M}$ is a finite free $S$-module equipped with filtration, Frobenius, monodromy, and an integrable connection $\nabla\colon \mathscr{M} \rightarrow \mathscr{M}\otimes_{\mathcal{O}_{L_0}} \widehat{\Omega}_{\mathcal{O}_{L_0} / W(k)}$ satisfying certain conditions (cf. Definition \ref{defn:rational-category-strly-div-lattices}). We naturally associate a lattice $T_{\mathrm{cris}}(\mathscr{M})$ of semistable $G_L$-representation with Hodge-Tate weights in $[0, r]$, and prove:

\begin{thm}[Theorem \ref{thm:equiv-category-strly-div-latt-semist-reps}] \label{thm:main-1}
The functor $T_{\mathrm{cris}}(\cdot)$ gives an anti-equivalence from the category of strongly divisible lattices of weight $r$ to the category of lattices of semistable $G_L$-representations with Hodge-Tate weights in $[0, r]$.
\end{thm}

We follow a similar strategy as in \cite{liu-semistable-lattice-breuil} to prove above theorem. The main new ingredient is the study of an integrable connection $\nabla\colon \mathscr{M} \rightarrow \mathscr{M}\otimes_{\mathcal{O}_{L_0}} \widehat{\Omega}_{\mathcal{O}_{L_0} / W(k)}$, which is crucial in understanding semistable representations of $G_L$. Such a connection is trivial in the case of perfect residue field. We adapt the idea in \cite{faltings} to study how this connection is precisely related to the Galois action. 

As an application, we consider in Section \ref{sec:cryst-cohom-strly-div-latt} a proper smooth formal scheme $\mathcal{X}$ over $\mathcal{O}_L$ and study its \'etale cohomology and crystalline cohomology. We give a cohomological description of the strongly divisible lattice associated to $H^i_{\text{\'et}}(\mathcal{X}_{\overline{L}}, \mathbf{Z}_p)$ for $i \leq p-2$, under the assumption that the crystalline cohomology of the special fiber of $\mathcal{X}$ is torsion-free in degrees $i$ and $i+1$. Denote the generic fiber of $\mathcal{X}$ by $X$, and let $\mathcal{X}_0 \coloneqq \mathcal{X}\times_{\mathcal{O}_L} \mathcal{O}_L/(p)$ and $\mathcal{X}_{k'} \coloneqq \mathcal{X}\times_{\mathcal{O}_L} k'$.

\begin{thm}[Theorem \ref{thm:cryst-cohom-strly-div-latt}] \label{thm-main-2}
Let $i \leq p-2$. Suppose that $H^i_{\mathrm{cris}}(\mathcal{X}_{k'} / \mathcal{O}_{L_0})$ and $H^{i+1}_{\mathrm{cris}}(\mathcal{X}_{k'} / \mathcal{O}_{L_0})$ are $p$-torsion free. Then

\begin{enumerate}
\item $T^i\coloneqq H^i_{\text{\'et}}(X_{\overline{L}}, \mathbf{Z}_p)$ is torsion free.
\item $\mathscr{M}^i \coloneqq H^i_{\mathrm{cris}}(\mathcal{X}_0 / S)$ is a strongly divisible lattice of weight $i$.
\item $T_{\mathrm{cris}}(\mathscr{M}^i) \cong (T^{i})^{\vee}$ as $G_L$-representations.	
\end{enumerate}	
\end{thm}
   
This generalizes a result of Cais-Liu \cite{cais-liu-breuil-kisin-mod-cryst-cohom}. We remark that the proof in \cite{cais-liu-breuil-kisin-mod-cryst-cohom} is based on $A_{\mathrm{inf}}$-cohomology theory in \cite{bhatt-morrow-scholze-integralpadic}. As we do not have yet the analogous $A_{\mathrm{inf}}$-cohomology theory when the residue field is imperfect, our argument is rather indirect and uses certain base change compatibilities. 
 
A main motivation for studying crystalline representations in the case of imperfect residue field (with a finite $p$-basis) is that it provides essential results for understanding crystalline local systems on any $p$-adic formal scheme with a smooth integral model over $\mathrm{Spf} \mathcal{O}_K$. For simplicity, consider the small affine case: let $R$ be the $p$-adic completion of an \'etale extension of $\mathcal{O}_K[T_1^{\pm 1}, \ldots, T_d^{\pm 1}]$. By the construction of de Rham and crystalline period rings in \cite{brinon-relative}, we have well-defined notions of de Rham \'etale local systems and crystalline \'etale local systems on the generic fiber of $\mathrm{Spf} R$. If we consider the $p$-adic completion $\widehat{R_{(\pi)}}$ of the localization $R_{(\pi)}$ and let $L_1 = \widehat{R_{(\pi)}}[p^{-1}]$, then $L_1$ is a complete discrete valuation field whose residue field has a finite $p$-basis. Let $\mathbb{L}$ be an isogeny $\mathbf{Z}_p$-local system on the generic fiber of $\mathrm{Spf} R$. Suppose $\mathbb{L}|_{\mathrm{Gal}(\overline{L_1}/L_1)}$ is crystalline. Then by \cite{tsuji-cryst-shvs} (see also \cite{moon-purity-crys-loc-syst}), $\mathbb{L}$ is crystalline. Thus, studying a crystalline \'etale local system on the generic fiber of $\mathrm{Spf} R$ amounts to studying an \'etale local system which is crystalline as a representation of ${\mathrm{Gal}(\overline{L_1}/L_1)}$. For example, some results established in Section \ref{sec:strongly-div-latt-semist-rep} in this article for crystalline case are used in proving the main theorem in \cite{Du-Liu-Moon-Shimizu-prism-F-cryst-loc-sys}.  

This paper is organized as follows. In Section \ref{sec:prelim-crys-rep}, we recall some preliminary results from \cite{brinon-crys-rep-imperfect-residue} and \cite{gao-integral-padic-hodge-imperfect} on crystalline representations and semistable representations in the imperfect residue field. In Section \ref{sec:etale-phi-mod}, we review the theory of \'etale $\varphi$-modules. In both Section \ref{sec:prelim-crys-rep} and Section \ref{sec:etale-phi-mod}, we also consider the compatibility of the constructions under certain base change maps, which will be used in the following sections. In Section \ref{sec:strongly-div-latt-semist-rep}, we define the category of strongly divisible lattices and construct a functor to lattices of semistable representations, and prove Theorem \ref{thm:main-1}. In Section \ref{sec:cryst-cohom-strly-div-latt}, we give a cohomological description of strongly divisible lattices for certain geometric situations and prove Theorem \ref{thm-main-2}.

\section*{Acknowledgements}

I would like to thank Bryden Cais and Tong Liu for many helpful discussions and communications during the preparation of this paper. I also thank the anonymous referee for valuable suggestions to improve the paper. The author was partially supported by AMS-Simons Travel Grant. Lastly, I would like to thank Soo Young Kim heartily for her continuous support.

\section{Crystalline and semistable representations in imperfect residue field case} \label{sec:prelim-crys-rep}

\subsection{Crystalline representations}

We follow the same notations as in Section \ref{sec:introduction}. As a preliminary, we review some results on crystalline $G_L$-representations over $\mathbf{Q}_p$ from \cite{brinon-crys-rep-imperfect-residue}. By a representation of $G_L$, we always mean a finite continuous representation. Choose a lifting $T_1, \ldots, T_d \in \mathcal{O}_{L_0}$ of a $p$-basis of $k'$. By \cite[Cor.~1.2.7(ii)]{berthelot-messing-dieudonne-cryst-III}, there exists a Frobenius endomorphism $\varphi$ on $\mathcal{O}_{L_0}$ such that $\varphi(T_i) = T_i^p$, and we fix such a Frobenius. Note that there exists a unique ring map $W(k) \rightarrow \mathcal{O}_{L_0}$ lifting $k \rightarrow k'$. By unicity, this map is compatible with Frobenius. Denote by $\widehat{\Omega}_{\mathcal{O}_{L_0}}$ the $p$-adically continuous Kahler differential of $\mathcal{O}_{L_0}$ relative to $\mathbf{Z}_p$, i.e. the $p$-adic completion of $\Omega^1_{\mathcal{O}_{L_0} / \mathbf{Z}_p}$. We have  
\[
\widehat{\Omega}_{\mathcal{O}_{L_0}} \cong \bigoplus_{i=1}^d \mathcal{O}_{L_0} \cdot d\log{T_i}.
\] 

We briefly recall crystalline period rings constructed in \cite{brinon-crys-rep-imperfect-residue}. Denote by $\mathcal{O}_{\overline{L}}^{\wedge}$ the $p$-adic completion of $\mathcal{O}_{\overline{L}}$, and let $\mathcal{O}_{\overline{L}}^{\flat} = \varprojlim_{\varphi} \mathcal{O}_{\overline{L}}/p\mathcal{O}_{\overline{L}}$ be its tilt. There exists a natural $W(k)$-linear surjective map $\theta\colon W(\mathcal{O}_{\overline{L}}^{\flat}) \rightarrow \mathcal{O}_{\overline{L}}^{\wedge}$ which lifts the projection onto the first factor. Let $\mathbf{A}_{\mathrm{cris}}(\mathcal{O}_{\overline{L}})$ be the $p$-adic completion of the PD-envelope of $W(\mathcal{O}_{\overline{L}}^{\flat})$ with respect to $\ker{\theta}$. $\mathbf{A}_{\mathrm{cris}}(\mathcal{O}_{\overline{L}})$ is equipped with $G_L$-action and Frobenius extending those on $W(\mathcal{O}_{\overline{L}}^{\flat})$, and is equipped with a natural PD-filtration. Let $\theta_{\mathcal{O}_{L_0}}\colon \mathcal{O}_{L_0}\otimes_{W(k)} W(\mathcal{O}_{\overline{L}}^{\flat}) \rightarrow \mathcal{O}_{\overline{L}}^{\wedge}$ be the $\mathcal{O}_{L_0}$-linear extension of $\theta$. The integral crystalline period ring $\mathbf{OA}_{\mathrm{cris}}(\mathcal{O}_{\overline{L}})$ is defined to be the $p$-adic completion of the PD-envelope of $\mathcal{O}_{L_0}\otimes_{W(k)} W(\mathcal{O}_{\overline{L}}^{\flat})$ with respect to $\ker{\theta_{\mathcal{O}_{L_0}}}$. $\mathbf{OA}_{\mathrm{cris}}(\mathcal{O}_{\overline{L}})$ is equipped with $G_L$-action, Frobenius, and PD-filtration compatible with those on $\mathbf{A}_{\mathrm{cris}}(\mathcal{O}_{\overline{L}})$. Moreover, we have a natural integrable conneciton
\[
\nabla\colon \mathbf{OA}_{\mathrm{cris}}(\mathcal{O}_{\overline{L}}) \rightarrow \mathbf{OA}_{\mathrm{cris}}(\mathcal{O}_{\overline{L}})\otimes_{\mathcal{O}_{L_0}} \widehat{\Omega}_{\mathcal{O}_{L_0}}
\]   
such that $\varphi$ is horizontal and $\mathbf{OA}_{\mathrm{cris}}(\mathcal{O}_{\overline{L}})^{\nabla = 0} = \mathbf{A}_{\mathrm{cris}}(\mathcal{O}_{\overline{L}})$. 

Choose compatibly $\epsilon_n \in \mathcal{O}_{\overline{L}}$ such that $\epsilon_0 =1, ~\epsilon_n = \epsilon_{n+1}^p$ with $\epsilon_1 \neq 1$, and let $\underline{\epsilon} = (\epsilon_n)_{n \geq 0} \in \mathcal{O}_{\overline{L}}^{\flat}$. Then $t \coloneqq \log{[\underline{\epsilon}]}\in A_{\mathrm{cris}}(\mathcal{O}_{\overline{L}})$. The crystalline period ring is defined to be $\mathbf{OB}_{\mathrm{cris}}(\mathcal{O}_{\overline{L}}) \coloneqq \mathbf{OA}_{\mathrm{cris}}(\mathcal{O}_{\overline{L}})[p^{-1}, t^{-1}]$. We also denote $\mathbf{B}_{\mathrm{cris}}(\mathcal{O}_{\overline{L}}) \coloneqq \mathbf{A}_{\mathrm{cris}}(\mathcal{O}_{\overline{L}})[p^{-1}, t^{-1}]$. For each $i = 1, \ldots, d$, we choose compatibly $T_{i, n} \in \mathcal{O}_{\overline{L}}$ such that $T_{i, 0} = T_i, ~T_{i, n} = T_{i, n+1}^p$, and write $\underline{T_i} = (T_{i, n})_{n \geq 0} \in \mathcal{O}_{\overline{L}}^{\flat}$. The following result on the structure of period rings is proved in \cite{brinon-crys-rep-imperfect-residue}.

\begin{lem}[cf.~{\cite[Prop.~2.39]{brinon-crys-rep-imperfect-residue}}] \label{lem:cryst-period-ring}
The map $X_i \mapsto T_i\otimes 1-1\otimes [\underline{T_i}]$ induces a $\mathbf{A}_{\mathrm{cris}}(\mathcal{O}_{\overline{L}})$-linear isomorphism
\[
\mathbf{A}_{\mathrm{cris}}(\mathcal{O}_{\overline{L}})\{X_1, \ldots, X_d\} \cong \mathbf{OA}_{\mathrm{cris}}(\mathcal{O}_{\overline{L}}),
\]	
where the former ring denotes the $p$-adically completed divided power polynomial with variables $X_i$ and coefficients in $\mathbf{A}_{\mathrm{cris}}(\mathcal{O}_{\overline{L}})$.   
\end{lem}

\begin{rem} \label{rem:period-ring}
We use different notations for crystalline period rings from \cite{brinon-crys-rep-imperfect-residue}, where $\mathbf{A}_{\mathrm{cris}}(\mathcal{O}_{\overline{L}})$, $\mathbf{OA}_{\mathrm{cris}}(\mathcal{O}_{\overline{L}})$, $\mathbf{B}_{\mathrm{cris}}(\mathcal{O}_{\overline{L}})$, $\mathbf{OB}_{\mathrm{cris}}(\mathcal{O}_{\overline{L}})$ are denoted by $A_{\mathrm{cris}}^{\nabla}$, $A_{\mathrm{cris}}$, $B_{\mathrm{cris}}^{\nabla}$, $B_{\mathrm{cris}}$ respectively.	
\end{rem}

For a $G_L$-representation $V$ over $\mathbf{Q}_p$, let 
\[
D_{\mathrm{cris}}(V) \coloneqq \mathrm{Hom}_{G_L}(V, \mathbf{OB}_{\mathrm{cris}}(\mathcal{O}_{\overline{L}})).
\]
Writing $V^{\vee} \coloneqq \mathrm{Hom}_{\mathbf{Q}_p}(V, \mathbf{Q}_p)$ for the dual representation of $V$, the natural morphism
\[
\alpha_{\mathrm{cris}}\colon D_{\mathrm{cris}}(V)\otimes_{L_0} \mathbf{OB}_{\mathrm{cris}}(\mathcal{O}_{\overline{L}})  \rightarrow V^{\vee}\otimes_{\mathbf{Q}_p}  \mathbf{OB}_{\mathrm{cris}}(\mathcal{O}_{\overline{L}})
\] 
is injective. We say $V$ is \textit{crystalline} if $\alpha_{\mathrm{cris}}$ is an isomorphism.

\begin{defn} \label{defn:fil-phi-nabla-mod}
A \textit{filtered} $(\varphi, \nabla)$-\textit{module} over $L_0$ is a finite $L_0$-module $D$ equipped with: 
\begin{itemize}
\item a $\varphi$-semi-linear injective emdomorphism $\varphi\colon D \rightarrow D$;
\item a decreasing filtration $\mathrm{Fil}^i D_L$ on $D_L \coloneqq D\otimes_{L_0} L$ by $L$-submodules such that $\mathrm{Fil}^i D_L = D_L$ for $i \ll 0$ and $\mathrm{Fil}^i D_L = 0$ for $i \gg 0$;
\item a topologically quasi-nilpotent integrable connection $\nabla\colon D \rightarrow D\otimes_{\mathcal{O}_{L_0}} \widehat{\Omega}_{\mathcal{O}_{L_0}}$ which satisfies Griffiths transversality and such that $\varphi$ is horizontal (cf. \cite[Def.~4.11, 4.4]{brinon-crys-rep-imperfect-residue}).	
\end{itemize}
Denote by $\mathrm{MF}(\varphi, \nabla)$ the category of filtered $(\varphi, \nabla)$-modules over $L_0$, whose morphisms are $L_0$-linear maps compatible with all structure.
\end{defn}

By \cite[Prop.~4.19]{brinon-crys-rep-imperfect-residue}, $D_{\mathrm{cris}}(\cdot)$ gives a functor from the category of $G_L$-representations over $\mathbf{Q}_p$ to $\mathrm{MF}(\varphi, \nabla)$. A filtered $(\varphi, \nabla)$-module is said to be \textit{weakly admissible} if it satisfies the usual conditions as in \cite[Def. 4.21]{brinon-crys-rep-imperfect-residue}. A main result in \cite{brinon-crys-rep-imperfect-residue} is the following:

\begin{thm} [cf.~{\cite[Cor.~4.37]{brinon-crys-rep-imperfect-residue}}] \label{thm:weakly-adm-cris}
The functor $D_{\mathrm{cris}}(\cdot)$ gives an anti-equivalence from the category of crystalline $G_L$-representations over $\mathbf{Q}_p$ to the category of weakly admissible filtered $(\varphi, \nabla)$-modules over $L_0$.	
\end{thm}

\begin{rem} \label{rem:contrav-Dcris}
We use the contravariant functor for $D_{\mathrm{cris}}(\cdot)$ whereas \cite{brinon-crys-rep-imperfect-residue} uses the covariant functor. A quasi-inverse of $D_{\mathrm{cris}}(\cdot)$ in the above theorem is given by
\[
V_{\mathrm{cris}}(D) \coloneqq \mathrm{Hom}_{\mathrm{Fil}, \varphi, \nabla}(D, \mathbf{OB}_{\mathrm{cris}}(\mathcal{O}_{\overline{L}})),
\]	
where the filtration compatibility means that the composite map 
\[
D_L \rightarrow \mathbf{OB}_{\mathrm{cris}}(\mathcal{O}_{\overline{L}})\otimes_{L_0} L \hookrightarrow \mathbf{OB}_{\mathrm{dR}}(\mathcal{O}_{\overline{L}})
\]
is compatible with filtrations. See \cite[Sec.~2.1]{brinon-crys-rep-imperfect-residue} for the construction of the de Rham period ring $\mathbf{OB}_{\mathrm{dR}}(\mathcal{O}_{\overline{L}})$ (denoted by $B_{\mathrm{dR}}$ in \textit{loc. cit.}), and also see  \cite[Prop.~2.47]{brinon-crys-rep-imperfect-residue}.
\end{rem}

For later use, we consider a certain base change map from $\mathcal{O}_L$ and functoriality of $D_{\mathrm{cris}}(\cdot)$. Let $R_{0, g}$ be the $p$-adic completion of the ring $\varinjlim_{\varphi}\mathcal{O}_{L_0}$, and let $k_g \coloneqq R_{0, g} / (p)$. The Frobenius on $\mathcal{O}_{L_0}$ extends uniquely to $R_{0, g}$, and we have a $\varphi$-compatible isomorphism $R_{0, g} \cong W(k_g)$. Let $b_g\colon \mathcal{O}_{L_0} \rightarrow W(k_g)$ be the induced map, and denote also by $b_g\colon L \rightarrow K_g \coloneqq W(k_g)\otimes_{W(k)} K$ the $K$-linear extension. Choose an embedding $\overline{L} \rightarrow \overline{K_g}$ extending $b_g$, which induces a map of Galois groups $G_{K_g} = \mathrm{Gal}(\overline{K_g}/K_g) \rightarrow G_L$. We then have an induced map of crystalline period rings $\mathbf{OB}_{\mathrm{cris}}(\mathcal{O}_{\overline{L}}) \rightarrow \mathbf{OB}_{\mathrm{cris}}(\mathcal{O}_{\overline{K_g}}) \cong \mathbf{B}_{\mathrm{cris}}(\mathcal{O}_{\overline{K_g}})$ compatible with all structures. In this way, if $V$ is a crystalline $G_L$-representation, then it can also be considered as a crystalline $G_{K_g}$-representation. Furthermore, for $D_{\mathrm{cris}, K_g}(V) \coloneqq \mathrm{Hom}_{G_{K_g}}(V, \mathbf{B}_{\mathrm{cris}}(\mathcal{O}_{\overline{K_g}}))$, we have 
\[
D_{\mathrm{cris}}(V)\otimes_{L_0, b_g} W(k_g)[p^{-1}] \cong D_{\mathrm{cris}, K_g}(V)
\]   
compatible with Frobenius and filtration.

\subsection{Semistable representations}

We have a natural notion of semistable representations of $G_L$ as in \cite{gao-integral-padic-hodge-imperfect}. Choose a uniformizer $\pi \in \mathcal{O}_K$. For integers $n \geq 0$, compatibly choose $\pi_n \in \overline{K}$ such that $\pi_0 = \pi$ and $\pi_{n+1}^p = \pi_n$, and let $\underline{\pi} = (\pi_n)_{n \geq 0} \in \mathcal{O}_{\overline{L}}^{\flat}$. Write $\mathfrak{u} \coloneqq \log[\underline{\pi}] \in \mathbf{B}_{\mathrm{dR}}(\mathcal{O}_{\overline{L}}) \cong \mathbf{B}_{\mathrm{dR}}(\mathcal{O}_{\overline{K_g}})$, and let $\mathbf{B}_{\mathrm{st}}(\mathcal{O}_{\overline{L}}) \coloneqq \mathbf{B}_{\mathrm{cris}}(\mathcal{O}_{\overline{L}})[\mathfrak{u}] \subset \mathbf{B}_{\mathrm{dR}}(\mathcal{O}_{\overline{L}})$. The semistable period ring is defined to be $\mathbf{OB}_{\mathrm{st}}(\mathcal{O}_{\overline{L}}) \coloneqq \mathbf{OB}_{\mathrm{cris}}(\mathcal{O}_{\overline{L}})[\mathfrak{u}]$, equipped with $\varphi$ given by $\varphi(\mathfrak{u}) = p\mathfrak{u}$ and $\mathbf{OB}_{\mathrm{cris}}(\mathcal{O}_{\overline{L}})$-linear derivation $N$ given by $N(\mathfrak{u}) = -1$.   

For a $G_L$-representation $V$ over $\mathbf{Q}_p$, let 
\[
D_{\mathrm{st}}(V) \coloneqq \mathrm{Hom}_{G_L}(V, \mathbf{OB}_{\mathrm{st}}(\mathcal{O}_{\overline{L}})).
\]
The natural morphism
\[
\alpha_{\mathrm{st}}\colon D_{\mathrm{st}}(V)\otimes_{L_0} \mathbf{OB}_{\mathrm{st}}(\mathcal{O}_{\overline{L}})  \rightarrow V^{\vee}\otimes_{\mathbf{Q}_p}  \mathbf{OB}_{\mathrm{st}}(\mathcal{O}_{\overline{L}})
\] 
is injective. We say $V$ is \textit{semistable} if $\alpha_{\mathrm{st}}$ is an isomorphism.

\begin{defn}
A \textit{filtered} $(\varphi, N, \nabla)$-\textit{module} over $L_0$ is a filtered $(\varphi, \nabla)$-\textit{module} $D$ equipped with $L_0$-linear map $N\colon D \rightarrow D$ which commutes with $\nabla$ and satisfies $N\varphi = p\varphi N$. Denote by $\mathrm{MF}(\varphi, N, \nabla)$ the category of filtered $(\varphi, N,  \nabla)$-modules over $L_0$ whose morphisms are $L_0$-linear maps compatible with all structure.	
\end{defn}

Note that $\mathrm{MF}(\varphi, \nabla)$ is a full subcategory of $\mathrm{MF}(\varphi, N, \nabla)$ consisting of objects with $N = 0$. Similarly as in \cite[Prop.~4.19 Pf.]{brinon-crys-rep-imperfect-residue}, $D_{\mathrm{st}}(\cdot)$ gives a functor from the category of $G_L$-representations over $\mathbf{Q}_p$ to $\mathrm{MF}(\varphi, N, \nabla)$. Weakly admissibility is defined as in the usual sense (cf.~\cite[Def.~2.4.8]{gao-integral-padic-hodge-imperfect}), and following is proved in \cite{gao-integral-padic-hodge-imperfect}. 

\begin{thm}[cf.~{\cite[Thm.~2.4.9]{gao-integral-padic-hodge-imperfect}}] \label{thm:weakly-adm-semist}
The functor $D_{\mathrm{st}}(\cdot)$ gives an anti-equivalence from the category of semistble $G_L$-representations over $\mathbf{Q}_p$ to the category of weakly admissible filtered $(\varphi, N, \nabla)$-modules over $L_0$. A quasi-inverse of $D_{\mathrm{st}}(\cdot)$ is given by
\[
V_{\mathrm{st}}(D) \coloneqq \mathrm{Hom}_{\mathrm{Fil}, \varphi, N, \nabla}(D, \mathbf{OB}_{\mathrm{st}}(\mathcal{O}_{\overline{L}})).
\]		
\end{thm}

As in the crystalline case, we have analogous functoriality for semistable representations under the base change $b_g\colon L \rightarrow K_g$.

\section{\'Etale $\varphi$-modules} \label{sec:etale-phi-mod}

We briefly recall some necessary facts on \'etale $\varphi$-modules from \cite[Sec.~7]{kim-groupscheme-relative} (cf. also \cite[Sec.~4.2]{gao-integral-padic-hodge-imperfect}). Let $\mathfrak{S} \coloneqq \mathcal{O}_{L_0}[\![u]\!]$ equipped with Frobenius given by $\varphi(u) = u^p$. Let $\mathcal{O}_{\mathcal{E}}$ be the $p$-adic completion of $\mathfrak{S}[u^{-1}]$, equipped with the Frobenius extending that on $\mathfrak{S}$.

\begin{defn} \label{defn:etale phi-module}
An \textit{\'{e}tale} $(\varphi, \mathcal{O}_\mathcal{E})$-\textit{module} is a pair $(\mathcal{M}, \varphi_{\mathcal{M}})$ where $\mathcal{M}$ is a finite free $\mathcal{O}_\mathcal{E}$-module and $\varphi_{\mathcal{M}}\colon \mathcal{M} \rightarrow \mathcal{M}$ is a $\varphi$-semilinear endomorphism such that $1\otimes\varphi_{\mathcal{M}}\colon \varphi^*\mathcal{M} \rightarrow \mathcal{M}$ is an isomorphism. Let $\mathrm{Mod}_{\mathcal{O}_{\mathcal{E}}}$ denote the category of \'{e}tale $(\varphi, \mathcal{O}_{\mathcal{E}})$-modules whose morphisms are $\mathcal{O}_\mathcal{E}$-linear maps compatible with Frobenius.
\end{defn}

We use \'{e}tale $(\varphi, \mathcal{O}_{\mathcal{E}})$-modules to study certain Galois representations. Let $L_{\infty} \coloneqq \bigcup_{n \geq 0} L(\pi_n) \subset \overline{L}$ and $\widetilde{L}_{\infty} \coloneqq \bigcup_{n \geq 0} L_{\infty}(T_{1, n}, \ldots, T_{d, n}) \subset \overline{L}$. Denote by $G_{L_\infty} \coloneqq \mathrm{Gal}(\overline{L}/L_{\infty})$ and $G_{\widetilde{L}_{\infty}} \coloneqq \mathrm{Gal}(\overline{L}/ \widetilde{L}_{\infty})$ the corresponding Galois subgroups of $G_L$. Let $\mathrm{Rep}_{\mathbf{Z}_p}(G_{\widetilde{L}_{\infty}})$ be the category of finite free $\mathbf{Z}_p$-representations of $G_{\widetilde{L}_{\infty}}$.

There exists a unique $W(k)$-linear map $\mathcal{O}_{L_0} \rightarrow W(\mathcal{O}_{\overline{L}}^{\flat})$ which maps $T_i$ to $[\underline{T_i}]$ and is compatible with Frobenius. We have a $\varphi$-equivariant embedding $\mathfrak{S} \rightarrow W(\mathcal{O}_{\overline{L}}^{\flat})$ given by $u \mapsto [\underline{\pi}]$, which extends to an embedding $\mathcal{O}_{\mathcal{E}} \rightarrow W(\overline{L}^{\flat})$. Let $\mathcal{O}_\mathcal{E}^{\mathrm{ur}}$ be the ring of integers of the maximal unramified extension of $\mathcal{O}_{\mathcal{E}}[p^{-1}]$ inside $W(\overline{L}^{\flat})[p^{-1}]$, and let $\widehat{\mathcal{O}}_{\mathcal{E}}^{\mathrm{ur}}$ be its $p$-adic completion. We define $\widehat{\mathfrak{S}}^{\mathrm{ur}} \coloneqq \widehat{\mathcal{O}}_{\mathcal{E}}^{\mathrm{ur}} \cap W(\mathcal{O}_{\overline{L}}^{\flat}) \subset W(\overline{L}^{\flat})$. The following is proved in \cite{kim-groupscheme-relative}.

\begin{lem} \label{lem:etale-ring-frob-galois} \emph{(cf.~\cite[Lem.~7.5, 7.6]{kim-groupscheme-relative})}
We have $(\widehat{\mathcal{O}}_{\mathcal{E}}^{\mathrm{ur}})^{G_{\widetilde{L}_\infty}} = \mathcal{O}_\mathcal{E}$. Furthermore, there exists a unique $G_{\widetilde{L}_\infty}$-equivariant ring endomorphism $\varphi$ on $\widehat{\mathcal{O}}_{\mathcal{E}}^{\mathrm{ur}}$ lifting the $p$-th power map on $\widehat{\mathcal{O}}_{\mathcal{E}}^{\mathrm{ur}}/(p)$ and extending $\varphi$ on $\mathcal{O}_\mathcal{E}$. The inclusion $\widehat{\mathcal{O}}_{\mathcal{E}}^{\mathrm{ur}} \hookrightarrow W(\overline{L}^{\flat})$ is $\varphi$-equivariant.
\end{lem}

For $\mathcal{M} \in \mathrm{Mod}_{\mathcal{O}_\mathcal{E}}$ and $T \in \mathrm{Rep}_{\mathbf{Z}_p}(G_{\widetilde{L}_\infty})$, define 
\[
T(\mathcal{M}) \coloneqq \mathrm{Hom}_{\mathcal{O}_\mathcal{E}, \varphi}(\mathcal{M}, \widehat{\mathcal{O}}_{\mathcal{E}}^{\mathrm{ur}}), ~~\mathcal{M}(T) \coloneqq \mathrm{Hom}_{G_{\widetilde{L}_\infty}}(T, \widehat{\mathcal{O}}_{\mathcal{E}}^{\mathrm{ur}}).
\]
The following equivalence is proved in \cite{kim-groupscheme-relative}.

\begin{prop} \label{prop:etale-gal-equiv} \emph{(cf.~\cite[Prop.~7.7]{kim-groupscheme-relative})} The assignments $T(\cdot)$ and $\mathcal{M}(\cdot)$ are exact rank-preserving anti-equivalences (quasi-inverse of each other) of $\otimes$-categories between $\mathrm{Mod}_{\mathcal{O}_{\mathcal{E}}}$ and $\mathrm{Rep}_{\mathbf{Z}_p}(G_{\widetilde{L}_\infty})$. Moreover, $T(\cdot)$ and $\mathcal{M}(\cdot)$ commute with taking duals. 
\end{prop}

We remark on the compatibility with respect to $b_g\colon L \rightarrow K_g$. Let $\mathfrak{S}_g \coloneqq W(k_g)[\![u]\!]$ and $\mathcal{O}_{\mathcal{E}, g}$ be the $p$-adic completion of $\mathfrak{S}_g[u^{-1}]$. Define $\widehat{\mathcal{O}}_{\mathcal{E}, g}^{\mathrm{ur}} \subset W(\overline{K_g}^{\flat})$ similarly as above and $\widehat{\mathfrak{S}}_g^{\mathrm{ur}} \coloneqq \widehat{\mathcal{O}}_{\mathcal{E}, g}^{\mathrm{ur}} \cap W(\mathcal{O}_{\overline{K_g}}^{\flat})$. Let $K_{g, \infty} \coloneqq \bigcup_{n \geq 0} K_g(\pi_n) \subset \overline{K_g}$, and write $G_{K_{g, \infty}} = \mathrm{Gal}(\overline{K_g} / K_{g, \infty})$. By \cite[Prop.~4.2.5 Pf.]{gao-integral-padic-hodge-imperfect}, under a suitable choice of embedding $\overline{L} \hookrightarrow \overline{K_g}$ extending $b_g$ so that $T_{i, n} \in K_g$ for each $i = 1, \ldots, d$ and $n \geq 0$, the map $G_{K_g} \rightarrow G_L$ induces an isomorphism $G_{K_{g, \infty}} \cong G_{\widetilde{L}_\infty}$. We fix such a choice from now on. Note that the $p$-adic completions of $\overline{L}$ and $\overline{K_g}$ are isomorphic, i.e. $\widehat{\overline{L}} = \widehat{\overline{K_g}}$, and thus $\overline{L}^{\flat} = \overline{K_g}^{\flat}$. If $\mathcal{M} \in \mathrm{Mod}_{\mathcal{O}_{\mathcal{E}}}$, then $\mathcal{M}_g = \mathcal{M}\otimes_{\mathcal{O}_{\mathcal{E}}} \mathcal{O}_{\mathcal{E}, g}$ with the induced tensor-product Frobenius $\varphi_{\mathcal{M}}\otimes \varphi_{\mathcal{O}_{\mathcal{E}, g}}$ is an \'etale $(\varphi, \mathcal{O}_{\mathcal{E}, g})$-module. Furthermore, we have a natural isomorphism
\[
T(\mathcal{M}) \cong T(\mathcal{M}_g) \coloneqq \mathrm{Hom}_{\mathcal{O}_{\mathcal{E}, g}, \varphi}(\mathcal{M}_g, \widehat{\mathcal{O}}_{\mathcal{E}, g}^{\mathrm{ur}})
\]  
of $G_{K_{g, \infty}}$-representations. 

\begin{lem} \label{lem:intersection-integral-closures}
We have	 
\[
\widehat{\mathcal{O}}_\mathcal{E}^{\mathrm{ur}} \cap \widehat{\mathfrak{S}}_g^{\mathrm{ur}} = \widehat{\mathfrak{S}}^{\mathrm{ur}}
\]
as subrings of $\widehat{\mathcal{O}}_{\mathcal{E}, g}^{\mathrm{ur}}$.
\end{lem}

\begin{proof}
Note that $\widehat{\mathcal{O}}_\mathcal{E}^{\mathrm{ur}} \subset W(\overline{L}^{\flat})$ and $\widehat{\mathfrak{S}}_g^{\mathrm{ur}} \subset W(\mathcal{O}_{\overline{K_g}}^{\flat}) = W(\mathcal{O}_{\overline{L}}^{\flat})$. So
\[
\widehat{\mathcal{O}}_\mathcal{E}^{\mathrm{ur}} \cap \widehat{\mathfrak{S}}_g^{\mathrm{ur}} = \widehat{\mathcal{O}}_\mathcal{E}^{\mathrm{ur}} \cap W(\mathcal{O}_{\overline{L}}^{\flat}) = \widehat{\mathfrak{S}}^{\mathrm{ur}}.
\] 	
\end{proof}

\section{Strongly divisible lattices and semistable representations} \label{sec:strongly-div-latt-semist-rep}

\subsection{Strongly divisible lattices} \label{sec:strongly-div-latt}

Fix a positive integer $r \leq p-2$. We will classify lattices of semistable $G_L$-representations with Hodge-Tate weights in $[0, r]$ via strongly divisible lattices, generalizing \cite[Thm.~2.3.5]{liu-semistable-lattice-breuil}. Denote by $\mathrm{MF}^{w, r}(\varphi, N, \nabla)$ the full subcategory of $\mathrm{MF}(\varphi, N, \nabla)$ consisting of weakly admissibly modules $D$ such that $\mathrm{Fil}^0 D_L = D_L$ and $\mathrm{Fil}^{r+1} D_L = 0$.

Let $E(u) \in W(k)[u]$ be the Eisenstein polynomial for $\pi$. Let $S$ be the $p$-adically completed divided power envelope of $\mathfrak{S}$ with respect to $(E(u))$. The Frobenius on $\mathfrak{S}$ extends uniquely to $S$. For each integer $i \geq 0$, let $\mathrm{Fil}^i S$ be the $p$-adically completed ideal of $S$ generated by the divided powers $\gamma_j(E(u)) = \frac{E(u)^j}{j!}, ~j \geq i$. Note that $\varphi(\mathrm{Fil}^i S) \subset p^i S$ if $i \leq p-1$. Let $N\colon S \rightarrow S$ be a $\mathcal{O}_{L_0}$-linear derivation given by $N(u) = -u$. We also have a natural integrable connection
\[
\nabla\colon S \rightarrow S\otimes_{\mathcal{O}_{L_0}} \widehat{\Omega}_{\mathcal{O}_{L_0}}
\]
given by the universal derivation on $\mathcal{O}_{L_0}$, which commutes with $N$. Note that the embedding  $\mathfrak{S} \rightarrow W(\mathcal{O}_{\overline{L}}^{\flat})$ extends to $S \rightarrow \mathbf{A}_{\mathrm{cris}}(\mathcal{O}_{\overline{L}})$, which is compatible with $\varphi$, filtration, and $G_{\widetilde{L}_{\infty}}$-action. Similarly as in \cite[Lem.~5.8 Pf.]{cais-liu-breuil-kisin-mod-cryst-cohom}, we see that the induced map $S/(p^n) \rightarrow \mathbf{A}_{\mathrm{cris}}(\mathcal{O}_{\overline{L}})/(p^n)$ is faithfully flat for each $n \geq 1$.

\begin{lem} \label{lem:filtration-S-Acris}
Let $x \in S$ \emph{(}resp. $\mathbf{A}_{\mathrm{cris}}(\mathcal{O}_{\overline{L}})$\emph{)}, and let $i, j \geq 0$ be any integers. If $E(u)^j x \in \mathrm{Fil}^{i+j} S$ \emph{(}resp. $E([\underline{\pi}])^j x \in \mathrm{Fil}^{i+j} \mathbf{A}_{\mathrm{cris}}(\mathcal{O}_{\overline{L}})$\emph{)}, then $x \in \mathrm{Fil}^{i} S$ \emph{(}resp. $x \in \mathrm{Fil}^{i} \mathbf{A}_{\mathrm{cris}}(\mathcal{O}_{\overline{L}})$\emph{)}.	
\end{lem}

\begin{proof}
This follows from a similar argument as in \cite[Lem.~3.2.2]{liu-semistable-lattice-breuil}.	
\end{proof}

For $D \in \mathrm{MF}^{w, r}(\varphi, N, \nabla)$, we define $\mathscr{D}(D) \coloneqq S[p^{-1}]\otimes_{L_0} D$ equipped with the tensor product Frobenius. Let $N\colon \mathscr{D}(D) \rightarrow \mathscr{D}(D)$ be the $L_0$-linear derivation given by $N_{S}\otimes 1+1\otimes N_D$, and let
\[
\nabla\colon \mathscr{D}(D) \rightarrow \mathscr{D}(D) \otimes_{\mathcal{O}_{L_0}} \widehat{\Omega}_{\mathcal{O}_{L_0}}
\] 
be the connection given by $\nabla_S\otimes 1+1\otimes \nabla_D$. For each $j = 1, \ldots, d$, let $N_{T_j}\colon \mathscr{D}(D) \rightarrow \mathscr{D}(D)$ be the derivation given by $\nabla\colon \mathscr{D}(D) \rightarrow \mathscr{D}(D) \otimes_{\mathcal{O}_{L_0}} \widehat{\Omega}_{\mathcal{O}_{L_0}} \cong \bigoplus_{j = 1}^d \mathscr{D}(D) \cdot d\log{T_j}$ composed with the projection to the $j$-th factor, and let $\partial_{T_j}\colon \mathscr{D}(D) \rightarrow \mathscr{D}(D)$ be the derivation given by $\partial_{T_j} = T_j^{-1}N_{T_j}$. Define a decreasing filtration on $\mathscr{D}(D)$ by $S[p^{-1}]$-submodules $\mathrm{Fil}^i \mathscr{D}(D)$ inductively by $\mathrm{Fil}^0 \mathscr{D}(D) = \mathscr{D}(D)$ and 
\[
\mathrm{Fil}^{i+1} \mathscr{D}(D) = \{x \in \mathscr{D}(D) ~|~ N(x) \in \mathrm{Fil}^i \mathscr{D}(D), ~q_{\pi}(x) \in \mathrm{Fil}^{i+1} D_L\}
\]
where $q_{\pi}\colon \mathscr{D}(D) \rightarrow D_L$ is the map induced by $S[p^{-1}] \rightarrow L, ~u \mapsto \pi$.

We thank Tong Liu for the following lemma on Griffiths transversality.

\begin{lem} \label{lem:Griffiths-transversality}
The connection $\nabla$ on $\mathscr{D}(D)$ satisfies Griffiths transversality:
\[
\nabla(\mathrm{Fil}^{i+1} \mathscr{D}(D)) \subset \mathrm{Fil}^{i} \mathscr{D}(D) \otimes_{\mathcal{O}_{L_0}} \widehat{\Omega}_{\mathcal{O}_{L_0}}.
\]	
\end{lem}

\begin{proof}
We need to show $N_{T_j}(\mathrm{Fil}^{i+1} \mathscr{D}(D)) \subset \mathrm{Fil}^{i} \mathscr{D}(D)$. We induct on $i$. The case $i = 0$ is clear. For $i \geq 1$, suppose $N_{T_j}(\mathrm{Fil}^{i} \mathscr{D}(D)) \subset \mathrm{Fil}^{i-1} \mathscr{D}(D)$, and let $x \in \mathrm{Fil}^{i+1} \mathscr{D}(D)$. We have
\[
N(N_{T_j}(x)) = N_{T_j}(N(x)) \in N_{T_j}(\mathrm{Fil}^{i} \mathscr{D}(D)) \subset \mathrm{Fil}^{i-1} \mathscr{D}(D). 
\]	
Furthermore, since $N_{T_j}\colon S \rightarrow S$ is $W(k)[\![u]\!]$-linear, we have $N_{T_j} \circ q_{\pi} = q_{\pi} \circ N_{T_j}$ as maps from $\mathscr{D}(D)$ to $D_L$. Thus,
\[
q_{\pi}(N_{T_j}(x)) = N_{T_j}(q_{\pi}(x)) \in N_{T_j}(\mathrm{Fil}^{i+1} D_L) \subset \mathrm{Fil}^{i} D_L.
\] 
Hence, $N_{T_j}(x) \in \mathrm{Fil}^{i} \mathscr{D}(D)$.
\end{proof}

\begin{defn} \label{defn:rational-category-strly-div-lattices}
Let $\mathrm{Mod}_{S}^r$ be a category whose objects are finite free $S$-modules $\mathscr{M}$ with $\mathscr{M}[p^{-1}] \cong \mathscr{D}(D)$ for some $D \in \mathrm{MF}^{w, r}(\varphi, N, \nabla)$, such that:
\begin{itemize}
\item $\mathscr{M}$ is stable under $\varphi_{\mathscr{D}(D)}$ and $N_{\mathscr{D}(D)}$;
\item $\mathscr{M}$ is stable under $\nabla_{\mathscr{D}(D)}$, i.e. $\nabla_{\mathscr{D}(D)}$ induces a connection $\nabla\colon \mathscr{M} \rightarrow \mathscr{M}\otimes_{\mathcal{O}_{L_0}} \widehat{\Omega}_{\mathcal{O}_{L_0}}$;
\item Let $\mathrm{Fil}^r \mathscr{M} \coloneqq \mathscr{M} \cap \mathrm{Fil}^r \mathscr{D}(D)$. Then $\varphi(\mathrm{Fil}^r \mathscr{M}) \subset p^r \mathscr{M}$. Denote $\varphi_r = \frac{\varphi}{p^r}\colon \mathrm{Fil}^r \mathscr{M} \rightarrow \mathscr{M}$.
\end{itemize}

Morphisms in $\mathrm{Mod}_{S}^r$ are $S$-linear maps compatible with $\varphi, N, \nabla, \mathrm{Fil}^r$. In above situation, we say $\mathscr{M}$ is a \textit{strongly divisible lattice of weight} $r$ in $\mathscr{D}(D)$. 	  	
\end{defn}

\begin{lem} \label{lem:phi_r}
Let $D \in \mathrm{MF}^{w, r}(\varphi, N, \nabla)$, and let $\mathscr{M} \in \mathrm{Mod}_{S}^r$ be a strongly divisible lattice of weight $r$ in $\mathscr{D}(D)$. Then $\mathscr{M}$ is generated by $\varphi_r(\mathrm{Fil}^r \mathscr{M})$ as $S$-modules. 	
\end{lem}

\begin{proof}
Since $D$ is weakly admissible, this follows from the same argument as in \cite[Prop. 2.1.3 Pf.]{breuil-semi-st-reps-strly-div-mod}.	
\end{proof}

We now construct a functor $T_{\mathrm{cris}}(\cdot)$ from $\mathrm{Mod}_{S}^r$ to the category of finite free $\mathbf{Z}_p$-representations of $G_L$. Let $\mathscr{M} \in \mathrm{Mod}_{S}^r$ and $D \in \mathrm{MF}^{w, r}(\varphi, N, \nabla)$ such that $\mathscr{M}$ is a strongly divisible lattice in $\mathscr{D}(D)$. We have a natural $\mathbf{A}_{\mathrm{cris}}(\mathcal{O}_{\overline{L}})$-semi-linear $G_{\widetilde{L}_{\infty}}$-action on $\mathbf{A}_{\mathrm{cris}}(\mathcal{O}_{\overline{L}})\otimes_{S} \mathscr{M}$ given by the $G_{\widetilde{L}_{\infty}}$-action on $\mathbf{A}_{\mathrm{cris}}(\mathcal{O}_{\overline{L}})$ and the trivial $G_{\widetilde{L}_{\infty}}$-action on $\mathscr{M}$. We extend this to a $G_L$-action using differential operators. For each $i = 1, \ldots, d$, let $N_{T_i}\colon \mathscr{M} \rightarrow \mathscr{M}$ be the derivation given by $\nabla\colon \mathscr{M} \rightarrow \mathscr{M}\otimes_{\mathcal{O}_{L_0}} \widehat{\Omega}_{\mathcal{O}_{L_0}}\cong \bigoplus_{i = 1}^d \mathscr{M} \cdot d\log{T_i}$ composed with the projection to the $i$-th factor. For any $\sigma \in G_L$, write
\[
\underline{\epsilon}(\sigma) \coloneqq \frac{\sigma([\underline{\pi}])}{[\underline{\pi}]} \text{ and }\underline{\mu_i}(\sigma) \coloneqq \frac{\sigma([\underline{T_i}])}{[\underline{T_i}]}, ~i = 1, \ldots, d.
\]
Note that $\log(\underline{\epsilon}(\sigma))$ and $\log(\underline{\mu_i}(\sigma))$ lie in $\mathrm{Fil}^1 \mathbf{A}_{\mathrm{cris}}(\mathcal{O}_{\overline{L}})$. For any $a\otimes x\in \mathbf{A}_{\mathrm{cris}}(\mathcal{O}_{\overline{L}})\otimes_{S} \mathscr{M}$, define
\begin{equation} \label{eq:Galois-action}
\sigma(a\otimes x) = \sum \sigma(a) \gamma_{i_0}(-\log(\underline{\epsilon}(\sigma)))\gamma_{i_1}(\log(\underline{\mu_1}(\sigma)))\cdots \gamma_{i_d}(\log(\underline{\mu_d}(\sigma))) \cdot N^{i_0}N_{T_1}^{i_1}\cdots N_{T_d}^{i_d}(x)
\end{equation}
where the sum goes over the multi-index $(i_0, i_1, \ldots, i_d)$ of non-negative integers and $\gamma_i$ denotes divided powers. Note that if $\sigma \in G_{\widetilde{L}_{\infty}}$, then $\sigma(a\otimes x) = \sigma(a)\otimes x$. Since $\nabla_{\mathscr{M}}$ and $N$ are topologically quasi-nilpotent and $\gamma_j(-\log(\underline{\epsilon}(\sigma))), ~\gamma_j(\log(\underline{\mu_i}(\sigma))) \rightarrow 0$ $p$-adically as $j \rightarrow \infty$, above sum converges. It is standard to check that this gives a well-defined $\mathbf{A}_{\mathrm{cris}}(\mathcal{O}_{\overline{L}})$-semi-linear $G_L$-action, which is compatible with $\varphi$ since $N\varphi = p\varphi N$ and $\varphi$ is horizontal with respect to $\nabla_{\mathscr{M}}$. Furthermore, this $G_L$-action is compatible with filtration since $N$ and $\nabla$ satisfy Griffiths transversality by definition and Lemma \ref{lem:Griffiths-transversality}.

Let
\[
T_{\mathrm{cris}}(\mathscr{M}) \coloneqq \mathrm{Hom}_{S, \mathrm{Fil}, \varphi} (\mathscr{M}, \mathbf{A}_{\mathrm{cris}}(\mathcal{O}_{\overline{L}})).
\]
Note that $T_{\mathrm{cris}}(\mathscr{M}) \cong \mathrm{Hom}_{\mathbf{A}_{\mathrm{cris}}(\mathcal{O}_{\overline{L}}), \mathrm{Fil}, \varphi} (\mathbf{A}_{\mathrm{cris}}(\mathcal{O}_{\overline{L}})\otimes_S \mathscr{M}, \mathbf{A}_{\mathrm{cris}}(\mathcal{O}_{\overline{L}}))$. Define the $G_L$-action on $T_{\mathrm{cris}}(\mathscr{M})$ by 
\[
\sigma(f)(x) = \sigma(f(\sigma^{-1}(x))) \text{ for any } x \in \mathbf{A}_{\mathrm{cris}}(\mathcal{O}_{\overline{L}})\otimes_S \mathscr{M}.
\]

To study $T_{\mathrm{cris}}(\mathscr{M})$, we consider $b_g\colon \mathcal{O}_{L_0} \rightarrow W(k_g)$. Let $S_g$ be the $p$-adically completed divided power envelope of $\mathfrak{S}_g$ with respect to $(E(u))$, equipped with $\varphi$, filtration, and $N$ similarly as above. Note that $D_g \coloneqq D\otimes_{L_0, b_g} W(k_g)[p^{-1}]$ is a weakly admissible filtered $(\varphi, N)$-module over $W(k_g)[p^{-1}]$, and $\mathscr{M}_g \coloneqq \mathscr{M}\otimes_{S, b_g} S_g$ with the induced tensor product $\varphi$, filtration, $N$ is a strongly divisible lattice of $S_g[p^{-1}]\otimes_{W(k_g)[p^{-1}]} D_g$ as in \cite{liu-semistable-lattice-breuil}. We have a natural injective map
\[
T_{\mathrm{cris}}(\mathscr{M}) \rightarrow T_{\mathrm{cris}}(\mathscr{M}_g) \coloneqq \mathrm{Hom}_{S_g, \mathrm{Fil}, \varphi}(\mathscr{M}_g, \mathbf{A}_{\mathrm{cris}}(\mathcal{O}_{\overline{K_g}})).
\]
Since the $G_{K_g}$-action on $T_{\mathrm{cris}}(\mathscr{M}_g)$ is given by \cite[Eq.~(5.1.1)]{liu-semistable-lattice-breuil}, this map is $G_{K_g}$-equivariant. Moreover, $T_{\mathrm{cris}}(\mathscr{M}_g)$ is finite free over $\mathbf{Z}_p$ by \cite{liu-semistable-lattice-breuil}, so $T_{\mathrm{cris}}(\mathscr{M})$ is finite free over $\mathbf{Z}_p$. We first study how this construction relates with \'etale $\varphi$-modules.

\begin{defn} \label{defn:quasi-kisin-mod}
A \textit{quasi-Kisin module} over $\mathfrak{S}$ of height $r$ is a finite free $\mathfrak{S}$-module $\mathfrak{M}$ equipped with $\varphi$-semi-linear endomorphism $\varphi\colon \mathfrak{M} \rightarrow \mathfrak{M}$ such that the cokernel of $1\otimes\varphi\colon \varphi^*\mathfrak{M} = \mathfrak{S}\otimes_{\varphi, \mathfrak{S}} \mathfrak{M} \rightarrow \mathfrak{M}$ is killed by $E(u)^r$.
\end{defn}

By Lemma \ref{lem:phi_r} and \cite[Thm.~5.1.3]{gao-integral-padic-hodge-imperfect}, there exists a quasi-Kisin module $\mathfrak{M}$ of height $r$ such that $S\otimes_{\varphi, \mathfrak{S}}\mathfrak{M} \cong \mathscr{M}$ compatible with Frobenius and
\[
\mathrm{Fil}^r \mathscr{M} = \{x \in S\otimes_{\varphi, \mathfrak{S}}\mathfrak{M} ~|~ (1\otimes\varphi_{\mathfrak{M}})(x) \in \mathrm{Fil}^r S \otimes_{\mathfrak{S}} \mathfrak{M}\}
\]
if we identify $S\otimes_{\varphi, \mathfrak{S}}\mathfrak{M} = \mathscr{M}$. Let $\mathcal{M} = \mathcal{O}_{\mathcal{E}}\otimes_{\mathfrak{S}}\mathfrak{M}$ with the induced tensor product Frobenius. Then $\mathcal{M}$ is an \'etale $\varphi$-module.

\begin{lem} \label{lem:etale-galois-isom}
The natural $G_{\widetilde{L}_\infty}$-equivariant map
\[
\mathrm{Hom}_{\mathfrak{S}, \varphi}(\mathfrak{M}, \widehat{\mathfrak{S}}^{\mathrm{ur}}) \rightarrow T(\mathcal{M}) = \mathrm{Hom}_{\mathcal{O}_\mathcal{E}, \varphi}(\mathcal{M}, \widehat{\mathcal{O}}_{\mathcal{E}}^{\mathrm{ur}})
\]	
is an isomorphism.
\end{lem}

\begin{proof}
By \cite[Sec.~B, Prop.~1.8.3]{fontaine-p-adic-rep-I}, the map
\[
\mathrm{Hom}_{\mathfrak{S}, \varphi}(\mathfrak{M}, \widehat{\mathfrak{S}}_{g}^{\mathrm{ur}}) \rightarrow \mathrm{Hom}_{\mathcal{O}_{\mathcal{E}}, \varphi}(\mathcal{M}, \widehat{\mathcal{O}}_{\mathcal{E}, g}^{\mathrm{ur}})
\]	
induced by $b_g\colon \mathcal{O}_{L_0} \rightarrow W(k_g)$ is an isomorphism. Since $T(\mathcal{M}) \cong \mathrm{Hom}_{\mathcal{O}_{\mathcal{E}}, \varphi}(\mathcal{M}, \widehat{\mathcal{O}}_{\mathcal{E}, g}^{\mathrm{ur}})$, the assertion follows from Lemma \ref{lem:intersection-integral-closures}.	
\end{proof}

Note that the embedding $\varphi\colon \widehat{\mathfrak{S}}^{\mathrm{ur}} \rightarrow \mathbf{A}_{\mathrm{cris}}(\mathcal{O}_{\overline{L}})$ induces a natural $G_{\widetilde{L}_\infty}$-equivariant injective map
\[
\mathrm{Hom}_{\mathfrak{S}, \varphi}(\mathfrak{M}, \widehat{\mathfrak{S}}^{\mathrm{ur}}) \rightarrow T_{\mathrm{cris}}(\mathscr{M}).
\]

\begin{lem} \label{lem:isom-galois-reps}
The natural maps $\mathrm{Hom}_{\mathfrak{S}, \varphi}(\mathfrak{M}, \widehat{\mathfrak{S}}^{\mathrm{ur}}) \rightarrow T_{\mathrm{cris}}(\mathscr{M})$ and $T_{\mathrm{cris}}(\mathscr{M}) \rightarrow T_{\mathrm{cris}}(\mathscr{M}_g)$ given above are $G_{\widetilde{L}_\infty}$-equivariant isomorphisms.	
\end{lem}

\begin{proof}
We have a commutative diagram
\[
\xymatrix{ 
\mathrm{Hom}_{\mathfrak{S}, \varphi}(\mathfrak{M}, \widehat{\mathfrak{S}}^{\mathrm{ur}}) \ar@{^{(}->}[r] \ar[d] &  T_{\mathrm{cris}}(\mathscr{M}) \ar@{^{(}->}[d]\\
\mathrm{Hom}_{\mathfrak{S}, \varphi}(\mathfrak{M}, \widehat{\mathfrak{S}}_{g}^{\mathrm{ur}}) \ar[r] & T_{\mathrm{cris}}(\mathscr{M}_g)
}
\]
The left vertical map is an isomorphism by Lemma \ref{lem:etale-galois-isom} and $T(\mathcal{M}) \cong T(\mathcal{M}_g)$, and the bottom horizontal map is an isomorphism by \cite[Lem.~3.3.4]{liu-semistable-lattice-breuil}. Since the other two maps are injective, they are isomorphisms.
\end{proof}  

In the remainder of this subsection, we show that $T_{\mathrm{cris}}(\cdot)$ gives a fully faithful functor from $\mathrm{Mod}_{S}^r$ to the category of lattices in semistable $G_L$-representations with Hodge-Tate weights in $[0, r]$. We first show that $T_{\mathrm{cris}}(\mathscr{M})[p^{-1}]$ is a semistable $G_L$-representation. By Theorem \ref{thm:weakly-adm-semist}, it suffices to show that 
\[
T_{\mathrm{cris}}(\mathscr{M})[p^{-1}] \cong V_{\mathrm{st}}(D) = \mathrm{Hom}_{\mathrm{Fil}, \varphi, N, \nabla}(D, \mathbf{OB}_{\mathrm{st}}(\mathcal{O}_{\overline{L}}))
\]
as $G_L$-representations. 

By \cite[Thm.~2.2.1, Prop.~3.5.1, Sec.~5.1]{liu-semistable-lattice-breuil}, we have a natural $G_{K_g}$-equivariant isomorphism 
\[
T_{\mathrm{cris}}(\mathscr{M}_g)[p^{-1}] \cong V_{\mathrm{st}}(D_g) = \mathrm{Hom}_{\mathrm{Fil}, \varphi, N}(D_g, \mathbf{B}_{\mathrm{st}}(\mathcal{O}_{\overline{K_g}})).
\]
Note that this isomorphism is induced by the natural maps
\begin{equation} \label{eq:map-Bst-hat}
D_g\otimes_{W(k_g)}\mathbf{B}_{\mathrm{st}}(\mathcal{O}_{\overline{K_g}}) \rightarrow D_g\otimes_{W(k_g)}\widehat{\mathbf{B}_{\mathrm{st}}} \cong \mathscr{M}_g\otimes_{S_g} \widehat{\mathbf{B}_{\mathrm{st}}} \rightarrow \mathscr{M}_g\otimes_{S_g} \mathbf{B}_{\mathrm{cris}}(\mathcal{O}_{\overline{K_g}}),
\end{equation}
which is $G_{K_g}$-equivariant by \cite[Lem.~5.2.1]{liu-semistable-lattice-breuil}. Here, $\widehat{\mathbf{B}_{\mathrm{st}}}$ denotes the period ring constructed in \cite{breuil-semi-st-reps-griffiths-trasv}, so that we have a natural embedding $\mathbf{B}_{\mathrm{st}}(\mathcal{O}_{\overline{K_g}}) \hookrightarrow \widehat{\mathbf{B}_{\mathrm{st}}}$ and projection $\widehat{\mathbf{B}_{\mathrm{st}}} \rightarrow \mathbf{B}_{\mathrm{cris}}(\mathcal{O}_{\overline{K_g}})$ which give above maps. Let $\iota_2\colon L_0 \rightarrow \mathbf{B}_{\mathrm{st}}(\mathcal{O}_{\overline{L}})$ be the map given by $T_i \mapsto [\underline{T_i}]$. Since $\mathbf{B}_{\mathrm{cris}}(\mathcal{O}_{\overline{L}}) = \mathbf{B}_{\mathrm{cris}}(\mathcal{O}_{\overline{K_g}})$ and $\mathscr{M}\otimes_S \mathbf{B}_{\mathrm{cris}}(\mathcal{O}_{\overline{L}}) \cong \mathscr{M}_g\otimes_{S_g} \mathbf{B}_{\mathrm{cris}}(\mathcal{O}_{\overline{K_g}})$, the composite of the map (\ref{eq:map-Bst-hat}) with $D\otimes_{L_0, \iota_2} \mathbf{B}_{\mathrm{st}}(\mathcal{O}_{\overline{L}}) \rightarrow D_g\otimes_{W(k_g)}\mathbf{B}_{\mathrm{st}}(\mathcal{O}_{\overline{K_g}})$ induces a $G_{K_g}$-equivariant isomorphism
\begin{equation} \label{eq:map-from-Tcris}
T_{\mathrm{cris}}(\mathscr{M})[p^{-1}] \cong \mathrm{Hom}_{\mathbf{B}_{\mathrm{st}}(\mathcal{O}_{\overline{L}}), \mathrm{Fil}, \varphi, N}(D\otimes_{L_0} \mathbf{B}_{\mathrm{st}}(\mathcal{O}_{\overline{L}}), \mathbf{B}_{\mathrm{st}}(\mathcal{O}_{\overline{L}})).
\end{equation}

We define $G_L$-action on $D\otimes_{L_0} \mathbf{B}_{\mathrm{st}}(\mathcal{O}_{\overline{L}})$ by 
\begin{equation} \label{eq:galois-action-D}
\sigma(x\otimes a) = \sum \sigma(a)\gamma_{i_1}(\log(\underline{\mu_1}(\sigma)))\cdots \gamma_{i_d}(\log(\underline{\mu_d}(\sigma))) \cdot N_{T_1}^{i_1}\cdots N_{T_d}^{i_d}(x)
\end{equation}
for $\sigma \in G_L$ and $x\otimes a \in D\otimes_{L_0} \mathbf{B}_{\mathrm{st}}(\mathcal{O}_{\overline{L}})$, where the sum goes over the multi-index $(i_1, \ldots, i_d)$ of non-negative integers. Then the map $D\otimes_{L_0} \mathbf{B}_{\mathrm{st}}(\mathcal{O}_{\overline{L}}) \rightarrow \mathscr{M}\otimes_S \mathbf{B}_{\mathrm{cris}}(\mathcal{O}_{\overline{L}})$ giving the isomorphism (\ref{eq:map-from-Tcris}) is also compatible with $G_{L_{\infty}}$-actions, since the $G_L$-action on $\mathscr{M}\otimes_S \mathbf{B}_{\mathrm{cris}}(\mathcal{O}_{\overline{L}})$ is given by Equation (\ref{eq:Galois-action}). Note that the Galois subgroups $G_{L_{\infty}}$ and $G_{K_g}$ generate $G_L$ by \cite[Lem.~4.4.1]{gao-integral-padic-hodge-imperfect}. Thus, the isomorphism (\ref{eq:map-from-Tcris}) is $G_L$-equivariant (where $\mathrm{Hom}_{\mathbf{B}_{\mathrm{st}}(\mathcal{O}_{\overline{L}}), \mathrm{Fil}, \varphi, N}(D\otimes_{L_0} \mathbf{B}_{\mathrm{st}}(\mathcal{O}_{\overline{L}}), \mathbf{B}_{\mathrm{st}}(\mathcal{O}_{\overline{L}}))$ is equipped with the $G_L$-action via $\sigma(f)(x) = \sigma(f(\sigma^{-1}(x)))$).  

\begin{lem} \label{lem:isom-Gal-rep-Bst-OBst}
We have a natural $G_L$-equivariant isomorphism
\[
\mathrm{Hom}_{\mathbf{B}_{\mathrm{st}}(\mathcal{O}_{\overline{L}}), \mathrm{Fil}, \varphi, N}(D\otimes_{L_0} \mathbf{B}_{\mathrm{st}}(\mathcal{O}_{\overline{L}}), \mathbf{B}_{\mathrm{st}}(\mathcal{O}_{\overline{L}})) \cong V_{\mathrm{st}}(D).
\]	
\end{lem}
 
\begin{proof}
By Lemma \ref{lem:cryst-period-ring}, we have a $\mathbf{B}_{\mathrm{st}}(\mathcal{O}_{\overline{L}})$-linear isomorphism
\[
\mathbf{B}_{\mathrm{st}}(\mathcal{O}_{\overline{L}})\{X_1, \ldots, X_d\} \cong \mathbf{OB}_{\mathrm{st}}(\mathcal{O}_{\overline{L}})
\]
given by $X_i \mapsto T_i\otimes 1-1\otimes [\underline{T_i}]$. Consider the projection
\[
pr\colon \mathbf{OB}_{\mathrm{st}}(\mathcal{O}_{\overline{L}}) \rightarrow \mathbf{B}_{\mathrm{st}}(\mathcal{O}_{\overline{L}})
\]
given by $X_i = T_i-[\underline{T_i}] \mapsto 0$, which is compatible with filtration, $\varphi$, and $N$. This induces the projection
\[
pr\colon D\otimes_{L_0, \iota_1} \mathbf{OB}_{\mathrm{st}}(\mathcal{O}_{\overline{L}}) \rightarrow D\otimes_{L_0, \iota_2} \mathbf{B}_{\mathrm{st}}(\mathcal{O}_{\overline{L}})
\] 
compatible with $\varphi, N$ and filtration (after tensoring with $L$ over $L_0$). Here, $\iota_1\colon L_0 \rightarrow  \mathbf{OB}_{\mathrm{st}}(\mathcal{O}_{\overline{L}})$ is the natural map given by $T_i \mapsto T_i\otimes 1$. We define a $\mathbf{B}_{\mathrm{st}}(\mathcal{O}_{\overline{L}})$-linear section $s$ to $pr$ as follows. For $x \in D$, let
\[
s(x) = \sum \gamma_{i_1}([\underline{T_1}]-T_1)\cdots \gamma_{i_d}([\underline{T_d}]-T_d)\cdot  \partial_{T_{1}}^{i_1}\cdots \partial_{T_d}^{i_d} (x)
\]
where the sum goes over the multi-index $(i_1, \ldots, i_d)$ of non-negative integers. Since $[\underline{T_i}]-T_i \in \mathrm{Fil}^1 \mathbf{OA}_{\mathrm{cris}}(\mathcal{O}_{\overline{L}})$ and $\partial_{T_i}$ is topologically quasi-nilpotent, this indeed defines a $\mathbf{B}_{\mathrm{st}}(\mathcal{O}_{\overline{L}})$-linear section $s\colon D\otimes_{L_0} \mathbf{B}_{\mathrm{st}}(\mathcal{O}_{\overline{L}}) \rightarrow D\otimes_{L_0} \mathbf{OB}_{\mathrm{st}}(\mathcal{O}_{\overline{L}})$ to $pr$. By \cite[Sec.~8.1]{Li-Liu-prismatic-cohom}, we have
\begin{equation} \label{eq:section}
s(x) = \sum (-1)^{i_1+\cdots+ i_d} \gamma_{i_1}(\log{(\frac{T_1}{[\underline{T_1}]})})\cdots \gamma_{i_d}(\log{(\frac{T_d}{[\underline{T_d}]})})\cdot N_{T_{1}}^{i_1}\cdots N_{T_d}^{i_d} (x).
\end{equation}
It follows from standard computations that $\nabla(s(x)) = 0$ and $\varphi(s(x)) = s(\varphi(x))$ for any $x \in D$. Moreover, $s$ is compatible with filtration since $\nabla_D$ satisfies Griffiths transversality. Note that $D$ has a $\mathcal{O}_{L_0}$-lattice stable under $\nabla_D$ on which the connection is topologically quasi-nilpotent. If $\{e_1, \ldots, e_m\}$ is a $\mathcal{O}_{L_0}$-basis of such a lattice, then $\{s(e_1), \ldots, s(e_m)\}$ generates $D\otimes_{L_0} \mathbf{OB}_{\mathrm{st}}(\mathcal{O}_{\overline{L}})$ as $\mathbf{OB}_{\mathrm{st}}(\mathcal{O}_{\overline{L}})$-modules. Thus, $s$ induces an isomorphism
\[
s\colon D\otimes_{L_0} \mathbf{B}_{\mathrm{st}}(\mathcal{O}_{\overline{L}}) \stackrel{\cong}{\rightarrow} (D\otimes_{L_0} \mathbf{OB}_{\mathrm{st}}(\mathcal{O}_{\overline{L}}))^{\nabla = 0}
\] 
of $\mathbf{B}_{\mathrm{st}}(\mathcal{O}_{\overline{L}})$-modules compatible with filtration and $\varphi, N$. 

By Equation (\ref{eq:galois-action-D}) and (\ref{eq:section}), we have
\[
pr(\sigma(s(x\otimes a))) = \sigma(x\otimes a) 
\]
for $\sigma \in G_L$ and $x\otimes a \in D\otimes_{L_0} \mathbf{B}_{\mathrm{st}}(\mathcal{O}_{\overline{L}})$. Thus, $s$ induces a $G_L$-equivariant isomorphism
\[
\mathrm{Hom}_{\mathbf{B}_{\mathrm{st}}(\mathcal{O}_{\overline{L}}), \mathrm{Fil}, \varphi, N}(D\otimes_{L_0} \mathbf{B}_{\mathrm{st}}(\mathcal{O}_{\overline{L}}), \mathbf{B}_{\mathrm{st}}(\mathcal{O}_{\overline{L}})) \cong V_{\mathrm{st}}(D).
\]
\end{proof}

Hence, from the isomorphism (\ref{eq:map-from-Tcris}) and Lemma \ref{lem:isom-Gal-rep-Bst-OBst}, we obtain $T_{\mathrm{cris}}(\mathscr{M})[p^{-1}] \cong V_{\mathrm{st}}(D)$ as $G_L$-representations. In particular, $T_{\mathrm{cris}}(\mathscr{M})$ is a lattice of the semistable representation $V_{\mathrm{st}}(D)$.

\begin{prop} \label{prop:Tcris-fully-faithful}
$T_{\mathrm{cris}}(\cdot)$ gives a fully faithful contravariant functor from $\mathrm{Mod}_{S}^r$ to the category of lattices in semistable $G_L$-representations with Hodge-Tate weights in $[0, r]$.
\end{prop}

\begin{proof}
It remains to show fully faithfulness. Let $\mathscr{M}_1, \mathscr{M}_2 \in \mathrm{Mod}_{S}^r$ and $D_1, D_2 \in \mathrm{MF}^{w, r}(\varphi, N, \nabla)$ such that $\mathscr{M}_i[p^{-1}] \cong \mathscr{D}(D_i), ~i = 1, 2$. By constructions above, we have a canonical isomorphism $T_{\mathrm{cris}}(\mathscr{M}_i)[p^{-1}] \cong V_{\mathrm{st}}(D_i)$ of $G_L$-representations. Identify $T_{\mathrm{cris}}(\mathscr{M}_i)$ as a lattice of $V_{\mathrm{st}}(D_i)$. 

Suppose we have a $G_L$-equivariant map $f\colon T_{\mathrm{cris}}(\mathscr{M}_1) \rightarrow T_{\mathrm{cris}}(\mathscr{M}_2)$. By Theorem \ref{thm:weakly-adm-semist}, we have a morphism $h\colon D_2 \rightarrow D_1$ in $\mathrm{MF}(\varphi, N, \nabla)$ such that $V_{\mathrm{st}}(h) = f[p^{-1}]$. Consider the induced map
\[
h\colon \mathscr{D}(D_2) = \mathscr{M}_2[p^{-1}] \rightarrow \mathscr{D}(D_1) = \mathscr{M}_1[p^{-1}].
\]
Note that $h$ is also compatible with filtration, since the induced map $\mathscr{D}(D_2)\otimes_{S, b_g} S_g \rightarrow \mathscr{D}(D_1)\otimes_{S, b_g} S_g$ is compatible with filtration by a main result in \cite{breuil-semi-st-reps-griffiths-trasv}. 
  
We need to show $h(\mathscr{M}_2) \subset \mathscr{M}_1$. By \cite[Cor.~3.5.2]{liu-semistable-lattice-breuil}, we have $h(\mathscr{M}_2) \subset \mathscr{M}_1 \otimes_{S} S_g$. So
\begin{align*}
h(\mathscr{M}_2) &\subset (\mathscr{M}_1\otimes_S S[p^{-1}]) \cap (\mathscr{M}_1 \otimes_{S} S_g) \\
	&= \mathscr{M}_1\otimes_S (S[p^{-1}] \cap S_g) = \mathscr{M}_1
\end{align*}
since $\mathscr{M}_1$ is free over $S$ and $S[p^{-1}] \cap S_g = S$. 	
\end{proof}

\subsection{Cartier Dual} \label{sec:cartier-dual}

To show the essential surjectivity of the functor $T_{\mathrm{cris}}(\cdot)$ constructed in Section \ref{sec:strongly-div-latt}, we first consider Cartier dual for strongly divisible lattices as in \cite[Sec.~4]{liu-semistable-lattice-breuil}. The results we need will follow essentially by the same arguments as in \textit{loc. cit}. It is convenient to study Cartier dual for objects in a category larger than $\mathrm{Mod}_S^{r}$. Let $c_1 = \frac{1}{p}\varphi(E(u)) \in S^{\times}$.

\begin{defn} \label{defn:category-quasi-strongly-divisible-lattices}
Let $\mathrm{Mod}_S^{r, q}$ be the category whose objects are tuples $(\mathscr{M}, \mathrm{Fil}^r \mathscr{M}, \varphi_r)$	such that
\begin{itemize}
\item $\mathscr{M}$ is a finite free $S$-module;
\item $\mathrm{Fil}^r \mathscr{M}$ is a $S$-submodule of $\mathscr{M}$ such that $\mathrm{Fil}^r S \cdot \mathscr{M} \subset \mathrm{Fil}^r \mathscr{M}$ and $\mathscr{M}/\mathrm{Fil}^r \mathscr{M}$ is $p$-torsion free;
\item $\varphi_r \colon \mathrm{Fil}^r \mathscr{M} \rightarrow \mathscr{M}$ is a $\varphi$-semi-linear map such that for any $s \in \mathrm{Fil}^r S$ and $x \in \mathscr{M}$, we have $\varphi_r(sx) = (c_1)^{-r}\varphi_r(s)\varphi_r(E(u)^r x)$;
\item $\varphi_r(\mathrm{Fil}^r \mathscr{M})$ generates $\mathscr{M}$ as $S$-modules.
\end{itemize}

Morphisms in $\mathrm{Mod}_S^{r, q}$ are $S$-linear maps compatible with all structures.
\end{defn}

\begin{lem} \label{lem:adapted-basis}
Let $\mathscr{M} \in \mathrm{Mod}_S^{r, q}$. There exist $\alpha_1, \ldots, \alpha_m \in \mathrm{Fil}^r \mathscr{M}$ such that
\begin{enumerate}
\item $\mathrm{Fil}^r \mathscr{M} = \bigoplus_{i=1}^m S\alpha_i+\mathrm{Fil}^p S\cdot \mathscr{M}$.
\item $E(u)^r \mathscr{M} \subset \bigoplus_{i=1}^m S\alpha_i$ and $\{\varphi_r(\alpha_1), \ldots, \varphi_{r}(\alpha_m)\}$ is a $S$-basis of $\mathscr{M}$.
\item If $A$ is the $m \times m$ matrix such that $(\alpha_1, \ldots, \alpha_m) = (\varphi_r(\alpha_1), \ldots, \varphi_{r}(\alpha_m)) A$, then there exists a $m\times m$ matrix $B$ with coefficients in $S$ such that $AB = E(u)^r I$.
\end{enumerate}
\end{lem}

\begin{proof}
This follows from the same argument as in \cite[Prop.~4.1.2 Pf.]{liu-semistable-lattice-breuil}, where the argument is valid without the assumption that the residue field of $L$ is perfect.	
\end{proof}

For $\mathscr{M} \in \mathrm{Mod}_S^{r, q}$, define its cartier dual $\mathscr{M}^* \coloneqq \mathrm{Hom}_S(\mathscr{M}, S), ~\mathrm{Fil}^r \mathscr{M}^* \coloneqq \{f \in \mathscr{M}^* ~|~ f(\mathrm{Fil}^r \mathscr{M}) \subset \mathrm{Fil}^r S\}$, and $\varphi_r\colon \mathrm{Fil}^r \mathscr{M}^* \rightarrow \mathscr{M}^*$ by $\varphi_r(f)(\varphi_r(x)) = \varphi_r(f(x))$ for any $x \in \mathrm{Fil}^r \mathscr{M}$. Note that $\varphi_r$ is well-defined since $\varphi_r(\mathrm{Fil}^r \mathscr{M}) = \mathscr{M}$	as $S$-modules. 

\begin{lem} \label{lem:cartier-dual-anti-equivalence}
The assignment $\mathscr{M} \rightarrow \mathscr{M}^*$ defines a functor from $\mathrm{Mod}_S^{r, q}$ to itself, which is an exact anti-equivalence such that $(\mathscr{M}^*)^* = \mathscr{M}$.
\end{lem}

\begin{proof}
This follows directly from the construction of the Cartier dual and Lemma \ref{lem:adapted-basis}.	
\end{proof}

\begin{eg} \label{eg:cartier-dual-of-S}
For the cartier dual $S^*$ of $S$, we have $\mathrm{Fil}^r S^* = S^*$ and $\varphi_r(1) = 1$. 
\end{eg}

For $\mathscr{M} \in \mathrm{Mod}_S^{r, q}$, let 
\[
T_{\mathrm{cris}}(\mathscr{M}) \coloneqq \mathrm{Hom}_{S, \mathrm{Fil}, \varphi_r} (\mathscr{M}, \mathbf{A}_{\mathrm{cris}}(\mathcal{O}_{\overline{L}})).
\]
As in Section \ref{sec:strongly-div-latt}, there exists a quasi-Kisin module $\mathfrak{M}$ of height $r$ such that $S\otimes_{\varphi, \mathfrak{S}}\mathfrak{M} \cong \mathscr{M}$ by \cite[Thm.~5.1.3]{gao-integral-padic-hodge-imperfect}, and the natural maps $\mathrm{Hom}_{\mathfrak{S}, \varphi}(\mathfrak{M}, \widehat{\mathfrak{S}}^{\mathrm{ur}}) \rightarrow T_{\mathrm{cris}}(\mathscr{M})$ and $T_{\mathrm{cris}}(\mathscr{M}) \rightarrow T_{\mathrm{cris}}(\mathscr{M}_g) \coloneqq \mathrm{Hom}_{S_g, \mathrm{Fil}, \varphi_r}(\mathscr{M}_g, \mathbf{A}_{\mathrm{cris}}(\mathcal{O}_{\overline{K_g}}))$ are $G_{\widetilde{L}_\infty}$-equivariant isomorphisms.

Consider the perfect pairing $\mathscr{M} \times \mathscr{M}^* \rightarrow S$ in the construction of Cartier dual, which is compatible with filtration and $\varphi$. Taking Cartier dual on both sides, we obtain a map
\[
S^* \rightarrow \mathscr{M}^* \times (\mathscr{M}^*)^* \cong \mathscr{M}^* \times \mathscr{M}
\]  
by Lemma \ref{lem:cartier-dual-anti-equivalence}. This induces a pairing
\[
T_{\mathrm{cris}}(\mathscr{M}) \times T_{\mathrm{cris}}(\mathscr{M}^*) \rightarrow \mathrm{Hom}_S (S^*, \mathbf{A}_{\mathrm{cris}}(\mathcal{O}_{\overline{L}})).
\]

\begin{lem} \label{lem:perfect-pairing}
The above pairing induces a perfect pairing of $\mathbf{Z}_p$-representations of $G_{\widetilde{L}_\infty}:$
\[
T_{\mathrm{cris}}(\mathscr{M}) \times T_{\mathrm{cris}}(\mathscr{M}^*) \rightarrow T_{\mathrm{cris}}(S^*) = \mathbf{Z}_p (r).
\]	
\end{lem}
  
\begin{proof}
It follows from the construction of Cartier dual that the image of the above pairing lies in $T_{\mathrm{cris}}(S^*)$. We have canonical $G_{\widetilde{L}_\infty}$-equivariant isomorphisms 
\[
T_{\mathrm{cris}}(\mathscr{M}) \cong T_{\mathrm{cris}}(\mathscr{M}_g), ~~ T_{\mathrm{cris}}(\mathscr{M}^*) \cong T_{\mathrm{cris}}(\mathscr{M}^*_g), ~~T_{\mathrm{cris}}(S^*) \cong T_{\mathrm{cris}}(S^*_g).
\] 	
Thus, the assertion follows by \cite[Lem.~4.3.1]{liu-semistable-lattice-breuil}.
\end{proof}

We will use the perfect pairing in above lemma to construct natural maps to relate $\mathscr{M}$ and $T_{\mathrm{cris}}(\mathscr{M})$. For $D \in \mathrm{MF}^{w, r}(\varphi, N, \nabla)$, define
\[
\mathrm{Fil}^r (\mathbf{A}_{\mathrm{cris}}(\mathcal{O}_{\overline{L}})\otimes_S \mathscr{D}(D)) \coloneqq \sum_{i = 0}^r \mathrm{Im}(\mathrm{Fil}^{r-i}\mathbf{A}_{\mathrm{cris}}(\mathcal{O}_{\overline{L}})\otimes_S \mathrm{Fil}^i \mathscr{D}(D)),
\]
where $\mathrm{Im}(\mathrm{Fil}^{r-i}\mathbf{A}_{\mathrm{cris}}(\mathcal{O}_{\overline{L}})\otimes_S \mathrm{Fil}^i \mathscr{D}(D))$ is the image of $\mathrm{Fil}^{r-i}\mathbf{A}_{\mathrm{cris}}(\mathcal{O}_{\overline{L}})\otimes_S \mathrm{Fil}^i \mathscr{D}(D)$ in $\mathbf{A}_{\mathrm{cris}}(\mathcal{O}_{\overline{L}})\otimes_S \mathscr{D}(D)$ under the natural map. If $\mathscr{M} \in \mathrm{Mod}_S^{r, q}$ such that $\mathscr{M}[p^{-1}] \cong \mathscr{D}(D)$ as $S[p^{-1}]$-modules compatibly with $\varphi$ and $\mathrm{Fil}^r \mathscr{M} = \mathscr{M} \cap \mathrm{Fil}^r \mathscr{D}(D)$, then define
\[
\mathrm{Fil}^r (\mathbf{A}_{\mathrm{cris}}(\mathcal{O}_{\overline{L}})\otimes_S \mathscr{M}) \coloneqq (\mathbf{A}_{\mathrm{cris}}(\mathcal{O}_{\overline{L}})\otimes_S \mathscr{M}) \cap \mathrm{Fil}^r (\mathbf{A}_{\mathrm{cris}}(\mathcal{O}_{\overline{L}})\otimes_S \mathscr{D}(D)).
\]  
Note that $\mathrm{Fil}^{r-i}S\cdot \mathrm{Fil}^i \mathscr{D}(D) \subset \mathrm{Fil}^r \mathscr{D}(D)$. Since the map $\mathrm{Fil}^r \mathscr{M}/p\mathrm{Fil}^r \mathscr{M} \rightarrow \mathscr{M}/p\mathscr{M}$ is injective and $S/(p) \rightarrow \mathbf{A}_{\mathrm{cris}}(\mathcal{O}_{\overline{L}})/(p)$ is flat, we deduce that $\mathrm{Fil}^r (\mathbf{A}_{\mathrm{cris}}(\mathcal{O}_{\overline{L}})\otimes_S \mathscr{M})$ is equal to the $p$-adic completion of $\mathrm{Im}(\mathbf{A}_{\mathrm{cris}}(\mathcal{O}_{\overline{L}})\otimes_S \mathrm{Fil}^r \mathscr{M})$. Let $\alpha_1, \ldots, \alpha_m \in \mathrm{Fil}^r \mathscr{M}$ as in Lemma \ref{lem:adapted-basis}. We then have
\begin{equation} \label{eq:filtration-r}
\mathrm{Fil}^r (\mathbf{A}_{\mathrm{cris}}(\mathcal{O}_{\overline{L}})\otimes_S \mathscr{M}) = \bigoplus_{i=1}^m \mathbf{A}_{\mathrm{cris}}(\mathcal{O}_{\overline{L}})\cdot \alpha_i + \mathrm{Fil}^p \mathbf{A}_{\mathrm{cris}}(\mathcal{O}_{\overline{L}}) \otimes_S \mathscr{M}.
\end{equation}
In particular, $\varphi_r\colon \mathrm{Fil}^r \mathscr{M} \rightarrow \mathscr{M}$ extends $\varphi_{\mathbf{A}_{\mathrm{cris}}(\mathcal{O}_{\overline{L}})}$-semi-linearly to
\[
\varphi_r\colon \mathrm{Fil}^r (\mathbf{A}_{\mathrm{cris}}(\mathcal{O}_{\overline{L}})\otimes_S \mathscr{M}) \rightarrow \mathbf{A}_{\mathrm{cris}}(\mathcal{O}_{\overline{L}})\otimes_S \mathscr{M}.
\]
 
Note that in above situation, $D^* \coloneqq D_{\mathrm{st}}(V_{\mathrm{st}}(D)^{\vee}(r))$ is in $\mathrm{MF}^{w, r}(\varphi, N, \nabla)$, and $\mathscr{M}^*[p^{-1}] \cong \mathscr{D}(D^*)$ compatibly with $\varphi$ and $\mathrm{Fil}^r$. We denote by $\mathbf{A}_{\mathrm{cris}}(\mathcal{O}_{\overline{L}})^*$ the ring $\mathbf{A}_{\mathrm{cris}}(\mathcal{O}_{\overline{L}})$ equipped with non-canonical filtration $\mathrm{Fil}^r \mathbf{A}_{\mathrm{cris}}(\mathcal{O}_{\overline{L}})^* = \mathbf{A}_{\mathrm{cris}}(\mathcal{O}_{\overline{L}})^*$ and Frobenius $\varphi_r(1) = 1$.

\begin{lem} \label{lem:Acris-Galois-reps}
We have natural $G_{\widetilde{L}_\infty}$-equivariant isomorphisms
\begin{align*}
\mathrm{Hom}_{\mathbf{A}_{\mathrm{cris}}(\mathcal{O}_{\overline{L}}), \mathrm{Fil}, \varphi} (\mathbf{A}_{\mathrm{cris}}(\mathcal{O}_{\overline{L}})^*, \mathbf{A}_{\mathrm{cris}}(\mathcal{O}_{\overline{L}})\otimes_S \mathscr{M}^*) &\cong \mathrm{Fil}^r(\mathbf{A}_{\mathrm{cris}}(\mathcal{O}_{\overline{L}})\otimes_S \mathscr{M}^*)^{\varphi_r = 1}\\
	&\cong \mathrm{Hom}_{S, \mathrm{Fil}, \varphi}(\mathscr{M}, \mathbf{A}_{\mathrm{cris}}(\mathcal{O}_{\overline{L}})).
\end{align*}	
\end{lem}

\begin{proof}
The first isomorphism is clear. For the second isomorphism, note that we have a natural isomorphism $\mathbf{A}_{\mathrm{cris}}(\mathcal{O}_{\overline{L}})\otimes_S \mathscr{M}^* \cong \mathrm{Hom}_{S}(\mathscr{M}, \mathbf{A}_{\mathrm{cris}}(\mathcal{O}_{\overline{L}}))$. Showing that this isomorphism induces the second isomorphism is based on Lemma \ref{lem:adapted-basis} and Equation (\ref{eq:filtration-r}), and it follows by the same argument as in \cite[Lem.~4.3.2 Pf.]{liu-semistable-lattice-breuil}. 	
\end{proof}

By above lemma, we obtain natural maps
\begin{equation} \label{eq:map-from-Tcris-to-Mdual}
T_{\mathrm{cris}}(\mathscr{M}) \cong \mathrm{Fil}^r(\mathbf{A}_{\mathrm{cris}}(\mathcal{O}_{\overline{L}})\otimes_S \mathscr{M}^*)^{\varphi_r = 1} \hookrightarrow \mathbf{A}_{\mathrm{cris}}(\mathcal{O}_{\overline{L}})\otimes_S \mathscr{M}^*,
\end{equation}  
and also a natural map $T_{\mathrm{cris}}(\mathscr{M}^*) \hookrightarrow \mathbf{A}_{\mathrm{cris}}(\mathcal{O}_{\overline{L}})\otimes_S \mathscr{M}$ similarly. Let $\mathfrak{t}$ be any $\mathbf{Z}_p$-basis of $\mathbf{Z}_p(1) = (\mathrm{Fil}^1 \mathbf{A}_{\mathrm{cris}}(\mathcal{O}_{\overline{L}}))^{\varphi_1 = 1} \subset \mathbf{A}_{\mathrm{cris}}(\mathcal{O}_{\overline{L}})$. We obtain a diagram
\begin{equation} \label{eq:commutative-diagram-Tcris-Mdual}
\xymatrix{
T_{\mathrm{cris}}(\mathscr{M}) \times T_{\mathrm{cris}}(\mathscr{M}^*) \ar@{^{(}->}[r] \ar[d] & \mathbf{A}_{\mathrm{cris}}(\mathcal{O}_{\overline{L}})\otimes_S \mathscr{M}^* \times \mathbf{A}_{\mathrm{cris}}(\mathcal{O}_{\overline{L}})\otimes_S \mathscr{M} \ar[d] \\
\mathbf{Z}_p(r) \ar[r]^{1 \mapsto \mathfrak{t}^r} & \mathbf{A}_{\mathrm{cris}}(\mathcal{O}_{\overline{L}})
}
\end{equation}
where the top horizontal map is given the map (\ref{eq:map-from-Tcris-to-Mdual}) and its dual, left vertical map by Lemma \ref{lem:perfect-pairing}, and right vertical map by the pairing $\mathscr{M}\times \mathscr{M}^* \rightarrow S$. Since the left vertical map is induced by taking dual of $\mathscr{M}\times \mathscr{M}^* \rightarrow S$, the above diagram is commutative.  

\begin{thm} \label{thm:comparing-M-Tcris}
There exist $\mathbf{A}_{\mathrm{cris}}(\mathcal{O}_{\overline{L}})$-linear injective maps
\begin{align*}
&\iota^*\colon T_{\mathrm{cris}}(\mathscr{M})^{\vee}(r)\otimes_{\mathbf{Z}_p} \mathbf{A}_{\mathrm{cris}}(\mathcal{O}_{\overline{L}})^* \hookrightarrow \mathbf{A}_{\mathrm{cris}}(\mathcal{O}_{\overline{L}})\otimes_S \mathscr{M},\\
&\iota\colon \mathbf{A}_{\mathrm{cris}}(\mathcal{O}_{\overline{L}})\otimes_S \mathscr{M} \hookrightarrow T_{\mathrm{cris}}(\mathscr{M})^{\vee}\otimes_{\mathbf{Z}_p} \mathbf{A}_{\mathrm{cris}}(\mathcal{O}_{\overline{L}})	
\end{align*}
which are compatible with $G_{\widetilde{L}_{\infty}}$-action, $\varphi$ and filtration, and satisfy $\iota\circ \iota^* = \mathrm{Id}\otimes \mathfrak{t}^r$.	
\end{thm}

\begin{proof}
Since $\mathbf{A}_{\mathrm{cris}}(\mathcal{O}_{\overline{L}}) = \mathbf{A}_{\mathrm{cris}}(\mathcal{O}_{\overline{K_g}})$, this follows from \cite[Thm.~4.3.4]{liu-semistable-lattice-breuil} for $\mathscr{M}_g = \mathscr{M}\otimes_S S_g$. For the convenience of reader, we sketch here the argument in \cite[Thm.~4.3.4 Pf.]{liu-semistable-lattice-breuil}.

Since $T_{\mathrm{cris}}(\mathscr{M}) = \mathrm{Hom}_{\mathbf{A}_{\mathrm{cris}}(\mathcal{O}_{\overline{L}}), \mathrm{Fil}, \varphi}(\mathbf{A}_{\mathrm{cris}}(\mathcal{O}_{\overline{L}})\otimes_S \mathscr{M}, \mathbf{A}_{\mathrm{cris}}(\mathcal{O}_{\overline{L}}))$, we have a natural map
\[
\tilde{\iota}\colon T_{\mathrm{cris}}(\mathscr{M}) \times \mathbf{A}_{\mathrm{cris}}(\mathcal{O}_{\overline{L}}) \times \mathscr{M} \rightarrow \mathbf{A}_{\mathrm{cris}}(\mathcal{O}_{\overline{L}}).	
\] 	
This induces a canonical map
\[
\iota\colon \mathbf{A}_{\mathrm{cris}}(\mathcal{O}_{\overline{L}})\otimes_S \mathscr{M} \rightarrow T_{\mathrm{cris}}(\mathscr{M})^{\vee}\otimes_{\mathbf{Z}_p} \mathbf{A}_{\mathrm{cris}}(\mathcal{O}_{\overline{L}})	
\]
which is compatible with $G_{\widetilde{L}_{\infty}}$-action, $\varphi$ and filtration. On the other hand, by Lemma \ref{lem:perfect-pairing} and the map (\ref{eq:map-from-Tcris-to-Mdual}), we get a canonical map
\[
\iota^*\colon T_{\mathrm{cris}}(\mathscr{M})^{\vee}(r)\otimes_{\mathbf{Z}_p} \mathbf{A}_{\mathrm{cris}}(\mathcal{O}_{\overline{L}})^* = T_{\mathrm{cris}}(\mathscr{M}^*)\otimes_{\mathbf{Z}_p} \mathbf{A}_{\mathrm{cris}}(\mathcal{O}_{\overline{L}})^* \rightarrow \mathbf{A}_{\mathrm{cris}}(\mathcal{O}_{\overline{L}})\otimes_S \mathscr{M}
\]
compatible with all structures. 

We obtain the following diagram
\[
\xymatrix{
T_{\mathrm{cris}}(\mathscr{M}) \times T_{\mathrm{cris}}(\mathscr{M}^*)\otimes_{\mathbf{Z}_p} \mathbf{A}_{\mathrm{cris}}(\mathcal{O}_{\overline{L}})^* \ar[r]^{~~~~~~\mathrm{Id}\otimes \iota^*} \ar[d] & T_{\mathrm{cris}}(\mathscr{M}) \times \mathbf{A}_{\mathrm{cris}}(\mathcal{O}_{\overline{L}})\otimes_S \mathscr{M} \ar[d]^{\tilde{\iota}} \\
\mathbf{Z}_p(r)\otimes_{\mathbf{Z}_p} \mathbf{A}_{\mathrm{cris}}(\mathcal{O}_{\overline{L}})^*  \ar[r]^{~~~~~~1 \mapsto \mathfrak{t}^r} & \mathbf{A}_{\mathrm{cris}}(\mathcal{O}_{\overline{L}})
}
\]
where the left vertical map is induced by Lemma \ref{lem:perfect-pairing}. The above diagram is commutative, since the map (\ref{eq:map-from-Tcris-to-Mdual}) is injective and the diagram (\ref{eq:commutative-diagram-Tcris-Mdual}) is commutative.
\end{proof}

\begin{lem} \label{lem:matrices-iota} 
Let $\alpha_1, \ldots, \alpha_m \in \mathrm{Fil}^r \mathscr{M}$ as in Lemma \ref{lem:adapted-basis}, and let $\{\mathfrak{e}_1, \ldots, \mathfrak{e}_m\}$ be a basis of $T_{\mathrm{cris}}(\mathscr{M})^{\vee}$. Let $C$ be the $m \times m$ matrix with coefficients in $\mathrm{Fil}^r \mathbf{A}_{\mathrm{cris}}(\mathcal{O}_{\overline{L}})$ such that $\iota(\alpha_1, \ldots, \alpha_m) = (\mathfrak{e}_1, \ldots, \mathfrak{e}_m) C$ 	given by Theorem \ref{thm:comparing-M-Tcris}. Then there exists a $m \times m$ matrix $C'$ with coefficients in $\mathbf{A}_{\mathrm{cris}}(\mathcal{O}_{\overline{L}})$ such that the coefficients of $C'C-\mathfrak{t}^r \mathrm{Id}$ are in $\mathrm{Fil}^p \mathbf{A}_{\mathrm{cris}}(\mathcal{O}_{\overline{L}})$.
\end{lem}

\begin{proof}
This follows from the same argument as in \cite[Lem.~4.3.6 Pf.]{liu-semistable-lattice-breuil}, based on Theorem \ref{thm:comparing-M-Tcris}, Equation (\ref{eq:filtration-r}), and Lemma \ref{lem:adapted-basis} and \ref{lem:filtration-S-Acris}.	
\end{proof}

\subsection{Essential Surjectivity of $T_{\mathrm{cris}}(\cdot)$} \label{sec:essential-surjectivity}

Let $V$ be a semistable $\mathbf{Q}_p$-representation of $G_L$ with Hodge-Tate weights in $[0, r]$, and let $T \subset V$ be a $G_L$-stable $\mathbf{Z}_p$-lattice. Denote $D = D_{\mathrm{st}}(V) \in \mathrm{MF}^{w, r}(\varphi, N, \nabla)$. By \cite[Thm.~3.2.3]{gao-integral-padic-hodge-imperfect}, there exists a quasi-Kisin module $\mathfrak{M}'$ over $\mathfrak{S}$ of height $r$ such that denoting $\varphi^* \mathfrak{M}' = \mathfrak{S}\otimes_{\varphi, \mathfrak{S}} \mathfrak{M}'$, we have $(\varphi^*\mathfrak{M}'/u\varphi^*\mathfrak{M}')[p^{-1}] \cong D$ compatibly with $\varphi$, and that $\mathfrak{M}'_g = \mathfrak{M}'\otimes_{\mathfrak{S}, b_g} \mathfrak{S}_g$ with the induced $\varphi$ is a Kisin module over $\mathfrak{S}_g$ associated with $D_g \coloneqq D\otimes_{L_0, b_g} W(k_g)[p^{-1}] = D_{\mathrm{st}}(V|_{G_{K_g}})$ via \cite[Lem.~1.3.13]{kisin-crystalline}. By \cite[Lem.~4.2.9]{gao-integral-padic-hodge-imperfect} and Lemma \ref{lem:etale-galois-isom}, 
\[
\mathfrak{M} \coloneqq \mathfrak{M}'[p^{-1}] \cap \mathcal{M}(T)
\]
is a quasi-Kisin module of height $r$ such that $\mathrm{Hom}_{\mathfrak{S}, \varphi}(\mathfrak{M}, \widehat{\mathfrak{S}}^{\mathrm{ur}}) \cong T$ as $G_{\widetilde{L}_{\infty}}$-representations. Let $\mathscr{M} \coloneqq S\otimes_{\varphi, \mathfrak{S}}\mathfrak{M}$ with the induced $\varphi$. 

\begin{rem}
If $\mathscr{M}_1 \in \mathrm{Mod}_{S}^r$ such that $T_{\mathrm{cris}}(\mathscr{M}_1) = T$ as $\mathbf{Z}_p$-representations of $G_L$, then the corresponding quasi-Kisin module of height $r$ given in Section \ref{sec:strongly-div-latt} agrees with $\mathfrak{M}$ above by \cite[Lem.~4.2.9]{gao-integral-padic-hodge-imperfect} and Lemma \ref{lem:isom-galois-reps}. 
\end{rem}

Consider the $\varphi$-compatible projection $q\colon S \twoheadrightarrow \mathcal{O}_{L_0}$ given by $u \mapsto 0$, and let $I_0 \subset S$ be its kernel and $M \coloneqq \mathscr{M}/I_0 \mathscr{M}$. We have the induced projection $q\colon \mathscr{M}[p^{-1}] \twoheadrightarrow M[p^{-1}] \cong D$.

\begin{lem} \label{lem:frob-compatible-section}
$q$ has a unique $\varphi$-compatible section $s\colon D \rightarrow \mathscr{M}[p^{-1}]$. Furthermore, $1\otimes s\colon S[p^{-1}]\otimes_{L_0} D \rightarrow \mathscr{M}[p^{-1}]$ is an isomorphism.	
\end{lem}

\begin{proof}
Since $\mathfrak{M}$ is a quasi-Kisin module of height $r$, the map $(1\otimes\varphi)[p^{-1}]\colon \varphi^*M[p^{-1}] = \mathcal{O}_{L_0}\otimes_{\varphi, \mathcal{O}_{L_0}}M[p^{-1}] \rightarrow M[p^{-1}]$ is an isomorphism, and the preimage of $M$ is contained in $p^{-r}(\varphi^* M)$. Thus, it follows from the standard argument as in \cite[Lem.~3.14 Pf.]{kim-groupscheme-relative} that there exists a unique $\varphi$-compatible section $s\colon D \cong M[p^{-1}] \rightarrow \mathscr{M}[p^{-1}]$. 

Consider the map $1\otimes s\colon S[p^{-1}]\otimes_{L_0} D \rightarrow \mathscr{M}[p^{-1}]$, which is a map of projective $S[p^{-1}]$-modules of the same rank. By Nakayama's lemma, there exists an element $a \in I_0 S[p^{-1}]$ such that the induced map
\[
S[p^{-1}]\otimes_{L_0}M[p^{-1}][(1+a)^{-1}] \rightarrow \mathscr{M}[p^{-1}][(1+a)^{-1}]
\]
is an isomorphism. On the other hand, by \cite[Lem.~1.2.6]{kisin-crystalline} (see also \cite[Sec.~3.2]{liu-semistable-lattice-breuil}), the map
\[
S_{g}[p^{-1}]\otimes_{L_0}M[p^{-1}] \rightarrow \mathscr{M}\otimes_S S_{g}[p^{-1}]
\]
induced by $b_g\colon \mathcal{O}_{L_0} \rightarrow W(k_g)$ is an isomorphism. We have
\[
S[p^{-1}][(1+a)^{-1}] \cap S_{g}[p^{-1}] = S[p^{-1}],
\]
since $L_0[\![u]\!] \cap S_{g}[p^{-1}] = S[p^{-1}]$ as subrings of $(W(k_g)[p^{-1}])[\![u]\!]$. Hence, the map $1\otimes s\colon S[p^{-1}]\otimes_{L_0}M[p^{-1}] \rightarrow \mathscr{M}[p^{-1}]$ is an isomorphism.
\end{proof}

By above lemma, we can identify $\mathscr{M}[p^{-1}] = \mathscr{D}(D)$ compatibly with $\varphi$. Define $\mathrm{Fil}^r \mathscr{M} = \mathscr{M} \cap \mathrm{Fil}^r \mathscr{D}(D)$. By \cite[Sec.~3.4]{liu-semistable-lattice-breuil}, $\mathscr{M}_g \coloneqq \mathscr{M}\otimes_S S_g$ is a quasi-strongly divisible lattice of $D_g\otimes_{W(k_g)} S_g[p^{-1}]$ as in \cite[Def.~2.3.3]{liu-semistable-lattice-breuil}. Since $S[p^{-1}] \cap S_g = S$ as submodules of $S_g[p^{-1}]$, we have
\[
\varphi_r(\mathrm{Fil}^r \mathscr{M}) \subset \mathscr{M}[p^{-1}] \cap \mathscr{M}_g = \mathscr{M}\otimes_S (S[p^{-1}] \cap S_g) = \mathscr{M},
\]
i.e. $\varphi(\mathrm{Fil}^r \mathscr{M}) \subset p^r \mathscr{M}$. Since $D$ is weakly admissible, we have $\varphi_r(\mathrm{Fil}^r \mathscr{M}) = \mathscr{M}$ as in Lemma \ref{lem:phi_r}. Thus, $\mathscr{M} \in \mathrm{Mod}_S^{r, q}$. 

Define the $G_L$-action on $\mathbf{A}_{\mathrm{cris}}(\mathcal{O}_{\overline{L}})\otimes_S \mathscr{D}(D)$ by the analogous formula as Equation (\ref{eq:Galois-action}), and let $V_{\mathrm{cris}}(\mathscr{D}(D))\coloneqq \mathrm{Hom}_{S, \mathrm{Fil}, \varphi}(\mathscr{D}(D), \mathbf{A}_{\mathrm{cris}}(\mathcal{O}_{\overline{L}})[p^{-1}])$ be equipped with the induced $G_L$-action. Similarly as in Section \ref{sec:strongly-div-latt} (by considering the corresponding maps induced via $b_g$), we have a natural $G_{\widetilde{L}_{\infty}}$-equivariant isomorphism
\[
T_{\mathrm{cris}}(\mathscr{M})[p^{-1}] \cong V_{\mathrm{cris}}(\mathscr{D}(D)).
\]
Furthermore, we have the natural $G_L$-equivariant isomorphism $V_{\mathrm{cris}}(\mathscr{D}(D)) \cong V$ by Section \ref{sec:strongly-div-latt}. Consider the commutative diagram
\begin{equation} \label{eq:commutative-diagram-rational-integral}
\xymatrix{
\mathbf{A}_{\mathrm{cris}}(\mathcal{O}_{\overline{L}})\otimes_S \mathscr{D}(D) \ar[r]^{\iota\otimes_{\mathbf{Z_p}}\mathbf{Q}_p~~~~~} & V_{\mathrm{cris}}(\mathscr{D}(D))^{\vee}\otimes_{\mathbf{Z}_p}\mathbf{A}_{\mathrm{cris}}(\mathcal{O}_{\overline{L}})\\
\mathbf{A}_{\mathrm{cris}}(\mathcal{O}_{\overline{L}})\otimes_S \mathscr{M} \ar[r]^{\iota~~~~~~} \ar@{^{(}->}[u] & T_{\mathrm{cris}}(\mathscr{M})^{\vee}\otimes_{\mathbf{Z}_p} \mathbf{A}_{\mathrm{cris}}(\mathcal{O}_{\overline{L}})	 \ar@{^{(}->}[u]
}
\end{equation}  
Note that the top horizontal map is $G_L$-equivariant, and the bottom map is $G_{\widetilde{L}_{\infty}}$-equivariant. Moreover, $T_{\mathrm{cris}}(\mathscr{M})^{\vee}\otimes_{\mathbf{Z}_p} \mathbf{A}_{\mathrm{cris}}(\mathcal{O}_{\overline{L}})$ is stable under $G_L$-action on $V_{\mathrm{cris}}(\mathscr{D}(D))^{\vee}\otimes_{\mathbf{Z}_p}\mathbf{A}_{\mathrm{cris}}(\mathcal{O}_{\overline{L}})$ since $T_{\mathrm{cris}}(\mathscr{M})^{\vee} = T^{\vee}$.

\begin{lem} \label{lem:p-divisibility-maps}
Let $\sigma$ be any element of $G_L$. Suppose we have a commutative diagram
\[
\xymatrix{
\mathbf{A}_{\mathrm{cris}}(\mathcal{O}_{\overline{L}})\otimes_S \mathscr{M} \ar[r]^{\iota~~~~~} \ar[d]^{f} & T_{\mathrm{cris}}(\mathscr{M})^{\vee}\otimes_{\mathbf{Z}_p} \mathbf{A}_{\mathrm{cris}}(\mathcal{O}_{\overline{L}}) \ar[d]^{f_1} \\
\mathbf{A}_{\mathrm{cris}}(\mathcal{O}_{\overline{L}})\otimes_S \mathscr{M} \ar[r]^{\iota~~~~~} & T_{\mathrm{cris}}(\mathscr{M})^{\vee}\otimes_{\mathbf{Z}_p} \mathbf{A}_{\mathrm{cris}}(\mathcal{O}_{\overline{L}}),
}
\] 	
where $f$ and $f_1$ are either $\mathbf{A}_{\mathrm{cris}}(\mathcal{O}_{\overline{L}})$-linear or $\sigma$-semi-linear morphisms compatible with $\varphi$ and filtration. If $p \mid f_1$, then $p \mid f$.
\end{lem}

\begin{proof}
Since $\mathbf{A}_{\mathrm{cris}}(\mathcal{O}_{\overline{L}}) = \mathbf{A}_{\mathrm{cris}}(\mathcal{O}_{\overline{K_g}})$, this follows from \cite[Lem.~5.3.1]{liu-semistable-lattice-breuil} for $\mathscr{M}_g = \mathscr{M}\otimes_S S_g$.
\end{proof}

\begin{cor} \label{cor:Galois-stable}
In the diagram \emph{(}\ref{eq:commutative-diagram-rational-integral}\emph{)}, $\mathbf{A}_{\mathrm{cris}}(\mathcal{O}_{\overline{L}})\otimes_S \mathscr{M}$ is stable under the $G_L$-action on $\mathbf{A}_{\mathrm{cris}}(\mathcal{O}_{\overline{L}})\otimes_S \mathscr{D}(D)$.	
\end{cor}

\begin{proof}
Let $\sigma$ be any element of $G_L$. Let $n \geq 0$ be an integer such that 
\[
p^n\sigma(\mathbf{A}_{\mathrm{cris}}(\mathcal{O}_{\overline{L}})\otimes_S \mathscr{M}) \subset \mathbf{A}_{\mathrm{cris}}(\mathcal{O}_{\overline{L}})\otimes_S \mathscr{M}.
\]
Since $T_{\mathrm{cris}}(\mathscr{M})^{\vee}\otimes_{\mathbf{Z}_p} \mathbf{A}_{\mathrm{cris}}(\mathcal{O}_{\overline{L}})$ is $G_L$-stable, we can apply Lemma \ref{lem:p-divisibility-maps} if $n \geq 1$ to deduce the assertion.	
\end{proof}

Now, to show that the functor $T_{\mathrm{cris}}(\cdot)$ in Section \ref{sec:strongly-div-latt} is essentially surjective, it suffices to show that $\mathscr{M}$ is stable under $N$ and $N_{T_j}$'s on $\mathscr{D}(D)$. By \cite[Lem.~3.5.3]{liu-semistable-lattice-breuil}, $N(\mathscr{M}_g) \subset \mathscr{M}_g$. So as above, we have
\[
N(\mathscr{M}) \subset \mathscr{D}(D) \cap \mathscr{M}_g = \mathscr{M}.
\]
It remains to show the stability under $N_{T_j}$ for each $j = 1, \ldots, d$. Recall that we chose an embedding $\overline{L} \hookrightarrow \overline{K_g}$ so that $G_{K_g}$ as a subgroup of $G_L$ acts trivially on $[\underline{T_j}]$ for each $j$. Denote $L_{p^{\infty}} = \bigcup_{n \geq 0} L(\epsilon_n)$ and $K_{g, p^{\infty}} = \bigcup_{n \geq 0} K_g(\epsilon_n)$. As explained in \cite[Not.~2.14]{brinon-crys-rep-imperfect-residue}, there exists $\tau_j \in \mathrm{Gal}(\overline{L}/L_{p^{\infty}})$ such that $\tau_j([\underline{T_j}]) = [\underline{\epsilon}][\underline{T_j}]$ (so $\log{(\underline{\mu_j}(\tau_j))}$ is a generator of $(\mathrm{Fil}^1 \mathbf{A}_{\mathrm{cris}}(\mathcal{O}_{\overline{L}}))^{\varphi_1 = 1}$) and that $\tau_j$ acts trivially on $[\underline{T_i}]$ if $i \neq j$. Since $p \geq 3$, by \cite[Lem.~5.1.2]{liu-semistable-lattice-breuil} applied to $G_{K_g}$, there exists $\sigma \in \mathrm{Gal}(\overline{K_g}/K_{g, p^{\infty}})$ such that $\sigma([\underline{\pi}]) = \tau_j([\underline{\pi}])$. Replacing $\tau_j$ by $\sigma^{-1}\tau_j$ if necessary, we may further assume that $\tau_j$ acts trivially on $[\underline{\pi}]$. Let $\mathfrak{t} = \log{(\underline{\mu_j}(\tau_j))}$.

Note that $\tau_j$ acts trivially on $\mathfrak{t}$. For any $x \in \mathscr{D}(D)$ and $n \geq 0$, we deduce by induction on $n$ that
\begin{equation} \label{eq:topological-generator}
(\tau_j-1)^n(x) = \sum_{m=n}^{\infty} (\sum_{i_1+\cdots+i_n = m, ~i_q \geq 1} \frac{m!}{i_1!\cdots i_n!}) \gamma_m(\mathfrak{t})\otimes N_{T_j}^m(x).
\end{equation}
We have $(\tau_j-1)^n(x) \in \mathrm{Fil}^n \mathbf{A}_{\mathrm{cris}}(\mathcal{O}_{\overline{L}})[p^{-1}]\otimes_S \mathscr{D}(D)$, and $\displaystyle \frac{1}{n}(\tau_j-1)^n(x) \rightarrow 0$ $p$-adically as $n \rightarrow \infty$. So 
\[
\log{(\tau_j)}(x) \coloneqq \sum_{n = 1}^\infty (-1)^{n-1}\frac{(\tau_j-1)^n (x)}{n} \in \mathbf{A}_{\mathrm{cris}}(\mathcal{O}_{\overline{L}})[p^{-1}]\otimes_S \mathscr{D}(D).
\]
By direct computation, we get
\[
\log{(\tau_j)}(x) = \mathfrak{t}\otimes N_{T_j}(x).
\]

We can now proceed similarly as in \cite[Sec.~5.3]{liu-semistable-lattice-breuil}. We first claim that $\mathfrak{t}\otimes N_{T_j}(\mathscr{M}) \subset \mathbf{A}_{\mathrm{cris}}(\mathcal{O}_{\overline{L}})\otimes_S \mathscr{M}$, for which it suffices to show $\displaystyle \frac{1}{n}(\tau_j-1)^n(\mathscr{M}) \subset \mathbf{A}_{\mathrm{cris}}(\mathcal{O}_{\overline{L}})\otimes_S \mathscr{M}$ for each $n \geq p$. Let $\alpha_1, \ldots, \alpha_m \in \mathrm{Fil}^r \mathscr{M}$ as in Lemma \ref{lem:adapted-basis} so that $\{\varphi_r(\alpha_1), \ldots, \varphi_r(\alpha_m)\}$ is a basis of $\mathscr{M}$. Since $\tau_j(\mathscr{M}) \subset \mathbf{A}_{\mathrm{cris}}(\mathcal{O}_{\overline{L}})\otimes_S \mathscr{M}$, Equation (\ref{eq:topological-generator}) gives
\[
(\tau_j-1)^n (\alpha_1, \ldots, \alpha_m) \in (\mathrm{Fil}^n \mathbf{A}_{\mathrm{cris}}(\mathcal{O}_{\overline{L}})\cdot (\mathbf{A}_{\mathrm{cris}}(\mathcal{O}_{\overline{L}})\otimes_S \mathscr{M}))^m.
\]
Thus,
\[
(\tau_j-1)^n (\varphi_r(\alpha_1), \ldots, \varphi_r(\alpha_m)) \in (\varphi_r(\mathrm{Fil}^n \mathbf{A}_{\mathrm{cris}}(\mathcal{O}_{\overline{L}}))\cdot \varphi(\mathbf{A}_{\mathrm{cris}}(\mathcal{O}_{\overline{L}})\otimes_S \mathscr{M}))^m.
\]
By \cite[Lem.~5.3.2]{liu-semistable-lattice-breuil}, $\displaystyle \frac{\varphi(E(u)^l)}{(np^r)l!} \in S$ for any integers $l\geq n \geq p$. So $\displaystyle \frac{1}{n}\varphi_r(\mathrm{Fil}^n \mathbf{A}_{\mathrm{cris}}(\mathcal{O}_{\overline{L}})) \subset \mathbf{A}_{\mathrm{cris}}(\mathcal{O}_{\overline{L}})$, which proves the claim.

Let $W$ be the $m \times m$ matrix with coefficients in $S[\frac{1}{p}]$ such that
\[
N_{T_j}(\varphi_r(\alpha_1), \ldots, \varphi_r(\alpha_m)) = (\varphi_r(\alpha_1), \ldots, \varphi_r(\alpha_m))W,
\]
and let $n\geq 0$ be the smallest integer such that the coefficients of $p^n W$ lie in $S$. Then $p^n N_{T_j}(\mathscr{M}) \subset \mathscr{M}$. Since $E(u)N_{T_j}(\mathrm{Fil}^r \mathscr{D}(D)) \subset \mathrm{Fil}^r \mathscr{D}(D)$ by Lemma \ref{lem:Griffiths-transversality}, we can write
\[
E(u)p^n N_{T_j}(\alpha_1, \ldots, \alpha_m) = (\alpha_1, \ldots, \alpha_m)X+(\varphi_r(\alpha_1), \ldots, \varphi_r(\alpha_m))Y
\]
for some $m \times m$ matrices $X$ and $Y$ with coefficients in $S$ and $\mathrm{Fil}^p S$ respectively by Lemma \ref{lem:adapted-basis}. On the other hand, since $\mathfrak{t}\otimes N_{T_j}(\mathscr{M}) \subset \mathbf{A}_{\mathrm{cris}}(\mathcal{O}_{\overline{L}})\otimes_S \mathscr{M}$, we have
\[
\mathfrak{t}\otimes N_{T_j}(\mathrm{Fil}^r \mathscr{M}) \subset (\mathbf{A}_{\mathrm{cris}}(\mathcal{O}_{\overline{L}})\otimes_S \mathscr{M}) \cap \mathrm{Fil}^r (\mathbf{A}_{\mathrm{cris}}(\mathcal{O}_{\overline{L}})\otimes_S \mathscr{D}(D)) = \mathrm{Fil}^r (\mathbf{A}_{\mathrm{cris}}(\mathcal{O}_{\overline{L}})\otimes_S \mathscr{M}). 
\]
So
\[
\mathfrak{t}N_{T_j}(\alpha_1, \ldots, \alpha_m) = (\alpha_1, \ldots, \alpha_m)X'+(\varphi_r(\alpha_1), \ldots, \varphi_r(\alpha_m))Y'
\]
for some $m \times m$ matrices $X'$ and $Y'$ with coefficients in $\mathbf{A}_{\mathrm{cris}}(\mathcal{O}_{\overline{L}})$ and $\mathrm{Fil}^p \mathbf{A}_{\mathrm{cris}}(\mathcal{O}_{\overline{L}})$ respectively. Let $A, B$ be $m \times m$ matrices as given in Lemma \ref{lem:adapted-basis}. Then we obtain
\[
E(u)^r(\mathfrak{t}X - E(u)p^n X') = -\mathfrak{t}BY + E(u)p^n BY'. 
\]
The matrix on the right hand side has coefficients in $\mathrm{Fil}^{p+1}\mathbf{A}_{\mathrm{cris}}(\mathcal{O}_{\overline{L}})$. So by Lemma \ref{lem:filtration-S-Acris}, $E(u)^{r-1}(\mathfrak{t}X - E(u)p^n X')$ has coefficients in $\mathrm{Fil}^{p}\mathbf{A}_{\mathrm{cris}}(\mathcal{O}_{\overline{L}})$. 

Suppose $n \geq 1$. Then the coefficients of $E(u)^{r-1}\mathfrak{t}X$ lie in $\mathrm{Fil}^{p}\mathbf{A}_{\mathrm{cris}}(\mathcal{O}_{\overline{L}})+p\mathbf{A}_{\mathrm{cris}}(\mathcal{O}_{\overline{L}})$. By a similar argument as in \cite[Lem.~5.3.1 Pf.]{liu-semistable-lattice-breuil} using Lemma \ref{lem:matrices-iota}, the coefficients of $X$ lie in $\mathrm{Fil}^{1}\mathbf{A}_{\mathrm{cris}}(\mathcal{O}_{\overline{L}})+p\mathbf{A}_{\mathrm{cris}}(\mathcal{O}_{\overline{L}})$. In particular, the coefficients of $\varphi(X)$ lie in $p\mathbf{A}_{\mathrm{cris}}(\mathcal{O}_{\overline{L}})$. Since $S/(p) \rightarrow \mathbf{A}_{\mathrm{cris}}(\mathcal{O}_{\overline{L}})/(p)$ is faithfully flat and so injective, the coefficients of $\varphi(X)$ lie in $pS$. Note that
\begin{align*}
c_1p^n N_{T_j}(\varphi_r(\alpha_1), \ldots, \varphi_r(\alpha_m)) &= p^n\varphi_r(E(u)N_{T_j}(\alpha_1, \ldots, \alpha_m))\\
	&= \varphi_r(\alpha_1, \ldots, \alpha_m)\varphi(X)+\varphi(\varphi_r(\alpha_1), \ldots, \varphi_r(\alpha_m))\varphi_r(Y)	.
\end{align*}
Since the coefficients of $\varphi(X)$ and $\varphi_r(Y)$ lie in $pS$, this contradicts the choice of $n$ to be minimal. Hence, $n = 0$.

This concludes the proof of the following theorem.

\begin{thm} \label{thm:equiv-category-strly-div-latt-semist-reps}
Let $0 \leq r \leq p-2$. The functor $T_{\mathrm{cris}}(\cdot)$ in Section \ref{sec:strongly-div-latt} gives an anti-equivalence from $\mathrm{Mod}_S^r$ to the category of $\mathbf{Z}_p$-lattices in semistable representations of $G_L$ with Hodge-Tate weights in $[0, r]$. 	
\end{thm}

Let $\mathrm{Mod}_S^{r, \mathrm{cris}}$ be the full subcategory of $\mathrm{Mod}_S^r$ consisting of $\mathscr{M}$ such that $\mathscr{M}[p^{-1}] \cong \mathscr{D}(D)$ for $D \in \mathrm{MF}^{w, r}(\varphi, N, \nabla)$ such that $N_D = 0$, satisfying the same conditions as in Definition \ref{defn:rational-category-strly-div-lattices}. Above arguments also imply the following for crystalline case.

\begin{thm} \label{thm:equiv-category-strly-div-latt-cris-reps}
Let $0 \leq r \leq p-2$. The functor $T_{\mathrm{cris}}(\cdot)$ in Section \ref{sec:strongly-div-latt} gives an anti-equivalence from $\mathrm{Mod}_S^{r, \mathrm{cris}}$ to the category of $\mathbf{Z}_p$-lattices in crystalline representations of $G_L$ with Hodge-Tate weights in $[0, r]$. 	
\end{thm}

\section{Crystalline cohomology and strongly divisible lattices} \label{sec:cryst-cohom-strly-div-latt}

Let $\mathcal{X}$ be a proper smooth formal scheme over $\mathcal{O}_L$. In this section, under certain assumptions, we study a cohomological description of the strongly divisible lattice associated with the \'etale cohomology $H^i_{\text{\'et}}(\mathcal{X}_{\overline{L}}, \mathbf{Z}_p)$ for $i \leq p-2$. More precisely, denote the generic fiber of $\mathcal{X}$ by $X$, and let $\mathcal{X}_0 \coloneqq \mathcal{X}\times_{\mathcal{O}_L} \mathcal{O}_L/(p)$ and $\mathcal{X}_{k'} \coloneqq \mathcal{X}\times_{\mathcal{O}_L} k'$. Let $T^i \coloneqq H^i_{\text{\'et}}(X_{\overline{L}}, \mathbf{Z}_p) / \mathrm{tors}$, and let $M^i\coloneqq H^i_{\mathrm{cris}}(\mathcal{X}_{k'} / \mathcal{O}_{L_0})$ and $\mathscr{M}^i \coloneqq H^i_{\mathrm{cris}}(\mathcal{X}_0 / S)$.

\begin{thm} \label{thm:cryst-cohom-strly-div-latt}
Let $i \leq p-2$. Suppose that $H^i_{\mathrm{cris}}(\mathcal{X}_{k'} / \mathcal{O}_{L_0})$ and $H^{i+1}_{\mathrm{cris}}(\mathcal{X}_{k'} / \mathcal{O}_{L_0})$ are $p$-torsion free. Then
\begin{enumerate}
\item $H^i_{\text{\'et}}(X_{\overline{L}}, \mathbf{Z}_p)$ is torsion free \emph{(}so $T^i = H^i_{\text{\'et}}(X_{\overline{L}}, \mathbf{Z}_p)$\emph{)}.
\item $\mathscr{M}^i \in \mathrm{Mod}_S^{i}$.
\item $T_{\mathrm{cris}}(\mathscr{M}^i) \cong (T^{i})^{\vee}$ as $G_L$-representations.	
\end{enumerate}	
\end{thm}

From now, we assume $M^i$ and $M^{i+1}$ are $p$-torsion free. To prove above theorem, we consider $b_g\colon L \rightarrow K_g$ as before. Since $b_g\colon \mathcal{O}_{L_0} \rightarrow W(k_g)$ is flat, we have $H^j_{\mathrm{cris}}(\mathcal{X}_{k'} / \mathcal{O}_{L_0})\otimes_{\mathcal{O}_{L_0}} W(k_g) \cong H^j_{\mathrm{cris}}(\mathcal{X}_{k_g} / W(k_g))$ for any $j$ by crystalline base change. So $H^i_{\mathrm{cris}}(\mathcal{X}_{k_g} / W(k_g))$ and $H^{i+1}_{\mathrm{cris}}(\mathcal{X}_{k_g} / W(k_g))$ are $p$-torsion free. In particular, we can apply \cite[Thm.~5.4]{cais-liu-breuil-kisin-mod-cryst-cohom} for $\mathcal{X}_{\mathcal{O}_{K_g}}$.

\begin{prop} \label{prop:etale-torsion-free}
$H^i_{\text{\'et}}(X_{\overline{L}}, \mathbf{Z}_p)$ is torsion free.
\end{prop}

\begin{proof}
Since $H^i_{\mathrm{cris}}(\mathcal{X}_{k_g} / W(k_g))$ is $p$-torsion free, we have $H^i_{\text{\'et}}(X_{\overline{K_g}}, \mathbf{Z}_p)$ is torsion free by \cite[Thm.~14.5]{bhatt-morrow-scholze-integralpadic}. On the other hand, by smooth and proper base change theorems, we have $H^i_{\text{\'et}}(X_{\overline{L}}, \mathbf{Z}_p) \cong H^i_{\text{\'et}}(X_{\overline{K_g}}, \mathbf{Z}_p)$. Thus, $H^i_{\text{\'et}}(X_{\overline{L}}, \mathbf{Z}_p)$ is torsion free. 
\end{proof}
 
Recall that $q\colon S \rightarrow \mathcal{O}_{L_0}$ is the projection given by $u \mapsto 0$. This induces a natural map $q\colon H^i_{\mathrm{cris}}(\mathcal{X}_0 / S) \rightarrow H^i_{\mathrm{cris}}(\mathcal{X}_{k'} / \mathcal{O}_{L_0})$.

\begin{prop} \label{prop:comparison-crys-cohom's}
There exists a unique section $s\colon H^i_{\mathrm{cris}}(\mathcal{X}_{k'} / \mathcal{O}_{L_0})[p^{-1}] \rightarrow H^i_{\mathrm{cris}}(\mathcal{X}_0 / S)[p^{-1}]$ of $q[p^{-1}]$ such that $s$ is $\varphi$-equivariant. Furthermore, the induced map 
\[
S\otimes_{\mathcal{O}_{L_0}} H^i_{\mathrm{cris}}(\mathcal{X}_{k'} / \mathcal{O}_{L_0})[p^{-1}] \rightarrow H^i_{\mathrm{cris}}(\mathcal{X}_0 / S)[p^{-1}]
\]
of $S[p^{-1}]$-modules is an isomorphism.	
\end{prop}

\begin{proof}
For each integer $n \geq 0$, set $u_n$ such that $u_0 = u$ and $u_{n+1}^p = u_n$. Recall that we choose $T_{i, n} \in \mathcal{O}_{\overline{L}}$ satisfying $T_{i, 0} = T_i, ~T_{i, n} = T_{i, n+1}^p$. Let $\mathcal{O}_{L_0, (n)} \coloneqq \mathcal{O}_{L_0}[T_{1, n}, \ldots, T_{d, n}]$ equipped with Frobenius given by $\varphi(T_{i, n}) = T_{i, n}^p$, which extends the Frobenius on $\mathcal{O}_{L_0}$. Equip $\mathcal{O}_{L_0, (n)}[u_n]$ with Frobenius given by $\varphi(u_n) = u_n^p$, and let $S_{(n)}$ be the $p$-adically completed PD-envelope of $\mathcal{O}_{L_0, (n)}[u_n]$ with respect to $(E(u_n))$. The Frobenius on $\mathcal{O}_{L_0, (n)}[u_n]$ extends naturally to $S_{(n)}$. Let $\mathcal{O}_{L, (n)} \coloneqq \mathcal{O}_{L_0, (n)}\otimes_{W(k)} \mathcal{O}_{K_n}$ where $K_n \coloneqq K[\pi_n]$. We have a natural inclusion $S \hookrightarrow S_{(n)}$. 

Let $e = [K: W(k)[p^{-1}]]$ be the ramification index. Consider the PD-thickenings $S_{(n)} \twoheadrightarrow \mathcal{O}_{L, (n)}/(\pi_n^e)$ given by $u_n \mapsto \pi_n$ and $S \twoheadrightarrow \mathcal{O}_L/(\pi^e) = \mathcal{O}_L/(p)$. Note that $\varphi^n\colon \mathcal{O}_{L_0, (n)} \rightarrow \mathcal{O}_{L_0}$ is an isomorphism, since $\{T_1, \ldots, T_d\}$ gives a $p$-basis of $\mathcal{O}_{L_0}/(p) = k'$ so the induced map $\varphi^n\colon \mathcal{O}_{L_0, (n)} / (p) \rightarrow \mathcal{O}_{L_0}/(p)$ is an isomorphism. Thus, $\varphi^n\colon S_{(n)} \rightarrow S$ is an isomorphism which is compatible with the isomorphism $\varphi^n\colon \mathcal{O}_{L, (n)}/(\pi_n^e) \stackrel{\cong}{\rightarrow} \mathcal{O}_L/(p)$. Denote $\mathcal{X}_{(n)} \coloneqq \mathcal{X}\times_{\mathcal{O}_L} \mathcal{O}_{L, (n)}/(\pi_n^e)$. By crystalline base change theorem, 
\begin{align*}
H^i_{\mathrm{cris}}(\mathcal{X}_{(n)} / S_{(n)})\otimes_{S_{(n)}, \varphi^n} S &\cong H^i_{\mathrm{cris}}(\mathcal{X}_{(n)}\times_{\mathcal{O}_{L, (n)}/(\pi_n^e), \varphi^n} \mathcal{O}_L/(p) / S)\\
	&\cong H^i_{\mathrm{cris}}(\mathcal{X}_0\times_{\mathcal{O}_L/(p), \varphi^n} \mathcal{O}_L/(p) / S).	
\end{align*}

Choose $n$ such that $p^n \geq e$. Then
\[
\mathcal{X}_{(n)} \cong \mathcal{X}_{\mathcal{O}_L/(\pi)}\times_{\mathcal{O}_L/(\pi)} \mathcal{O}_{L, (n)}/(\pi_n^e) = \mathcal{X}_{k'}\times_{k'} \mathcal{O}_{L, (n)}/(\pi_n^e).
\]	
The natural inclusion $\mathcal{O}_{L_0} \hookrightarrow S_{(n)}$ is a PD-morphism over $\mathcal{O}_{L_0}/(p) \rightarrow \mathcal{O}_{L, (n)}/(\pi_n^e)$. Since $\mathcal{O}_{L_0} \hookrightarrow S_{(n)}$ is flat, we have by crystalline base change theorem that 
\[
H^i_{\mathrm{cris}}(\mathcal{X}_{k'} / \mathcal{O}_{L_0})\otimes_{\mathcal{O}_{L_0}} S_{(n)} \cong H^i_{\mathrm{cris}}(\mathcal{X}_{(n)}/S_{(n)}).
\]

From the above isomorphisms,
\[
H^i_{\mathrm{cris}}(\mathcal{X}_{k'} / \mathcal{O}_{L_0})\otimes_{\mathcal{O}_{L_0}, \varphi^n} S \cong H^i_{\mathrm{cris}}(\mathcal{X}_0\times_{\mathcal{O}_L/(p), \varphi^n} \mathcal{O}_L/(p) / S).	
\]
Thus, we obtain a morphism
\begin{equation} \label{eq:map-cryst-cohom}
H^i_{\mathrm{cris}}(\mathcal{X}_{k'} / \mathcal{O}_{L_0})\otimes_{\mathcal{O}_{L_0}, \varphi^n} \mathcal{O}_{L_0} \rightarrow H^i_{\mathrm{cris}}(\mathcal{X}_0\times_{\mathcal{O}_L/(p), \varphi^n} \mathcal{O}_L/(p) / S) \rightarrow H^i_{\mathrm{cris}}(\mathcal{X}_0/S)	
\end{equation}
where the second map is the $n$-th iteration of relative Frobenius. Note that by \cite[Thm.~1.1]{bhatt-morrow-scholze-integralpadic}, $H^i_{\mathrm{cris}}(\mathcal{X}_{k_g} / W(k_g))[p^{-1}]$ is a filtered $\varphi$-module over $W(k_g)[p^{-1}]$, and so $\varphi^n\colon H^i_{\mathrm{cris}}(\mathcal{X}_{k_g} / W(k_g))[p^{-1}] \rightarrow H^i_{\mathrm{cris}}(\mathcal{X}_{k_g} / W(k_g))[p^{-1}]$ is injective. Since $H^i_{\mathrm{cris}}(\mathcal{X}_{k'} / \mathcal{O}_{L_0})\otimes_{\mathcal{O}_{L_0}} W(k_g) \cong H^i_{\mathrm{cris}}(\mathcal{X}_{k_g} / W(k_g))$, we have that $\varphi^n\colon H^i_{\mathrm{cris}}(\mathcal{X}_{k'} / \mathcal{O}_{L_0})[p^{-1}] \rightarrow H^i_{\mathrm{cris}}(\mathcal{X}_{k'} / \mathcal{O}_{L_0})[p^{-1}]$ is injective. Thus, the map 
\[
\varphi^n\otimes 1\colon H^i_{\mathrm{cris}}(\mathcal{X}_{k'} / \mathcal{O}_{L_0})[p^{-1}]\otimes_{\mathcal{O}_{L_0}, \varphi^n} \mathcal{O}_{L_0} \rightarrow H^i_{\mathrm{cris}}(\mathcal{X}_{k'} / \mathcal{O}_{L_0})[p^{-1}]
\]
is an isomorphism. Composing the inverse of this isomorphism with the map (\ref{eq:map-cryst-cohom}) with $p$ inverted, we obtain a map
\[
s\colon H^i_{\mathrm{cris}}(\mathcal{X}_{k'} / \mathcal{O}_{L_0})[p^{-1}] \rightarrow H^i_{\mathrm{cris}}(\mathcal{X}_0 / S)[p^{-1}].
\]

It is clear from the constructions that $s$ is a $\varphi$-equivariant section of $q[p^{-1}]$. The uniqueness follows from the standard argument. By \cite[Thm.~1.6]{berthelot-ogus-F-isocrys-deRham-cohom}, the second map in (\ref{eq:map-cryst-cohom}) with $p$ inverted is an isomorphism. Hence, the induced map  
\[
S\otimes_{\mathcal{O}_{L_0}} H^i_{\mathrm{cris}}(\mathcal{X}_{k'} / \mathcal{O}_{L_0})[p^{-1}]  \rightarrow H^i_{\mathrm{cris}}(\mathcal{X}_0 / S)[p^{-1}]
\]
of $S[p^{-1}]$-modules is an isomorphism.
\end{proof}

We now consider the natural connection on crystalline cohomology. Let $\mathcal{O}_{L_0}(1)$ be the $p$-adically completed PD-envelope of $\mathcal{O}_{L_0}\widehat{\otimes}_{W(k)}\mathcal{O}_{L_0}$ with respect to the kernel of the natural map $\mathcal{O}_{L_0}\widehat{\otimes}_{W(k)}\mathcal{O}_{L_0} \twoheadrightarrow \mathcal{O}_{L_0}/(p)$, and let $S(1)$ be the $p$-adically completed PD-envelope of $\mathcal{O}_{L_0}\widehat{\otimes}_{W(k)}\mathfrak{S}$ with respect to the kernel of the map $\mathcal{O}_{L_0}\widehat{\otimes}_{W(k)}\mathfrak{S} \twoheadrightarrow \mathcal{O}_{L}/(p)$ given by $u \mapsto \pi$. Here, $\widehat{\otimes}$ denotes the $p$-adically completed $\otimes$-product. Note that the kernel of $\mathcal{O}_{L_0}\widehat{\otimes}_{W(k)}\mathcal{O}_{L_0} \twoheadrightarrow \mathcal{O}_{L_0}/(p)$ is generated by $p$ and $a\otimes 1 - 1\otimes a, ~a \in \mathcal{O}_{L_0}$, and the kernel of $\mathcal{O}_{L_0}\widehat{\otimes}_{W(k)}\mathfrak{S} \twoheadrightarrow \mathcal{O}_{L}/(p)$ is generated by $p$, $E(u)$, and $a\otimes 1 - 1\otimes a, ~a \in \mathcal{O}_{L_0}$. The map $q\colon S \rightarrow \mathcal{O}_{L_0}$ naturally extends to $S(1) \rightarrow \mathcal{O}_{L_0}(1)$ via $u \mapsto 0$.

Let $p_1, p_2\colon \mathcal{O}_{L_0} \rightarrow \mathcal{O}_{L_0}(1)$ be the maps given by
\[
p_1(a) = a\otimes 1, ~~p_2(a) = 1\otimes a, ~~a \in \mathcal{O}_{L_0}.
\]
We also denote by $p_1, p_2\colon \mathfrak{S} \rightarrow S(1)$ the maps given by
\[
p_1(\sum_{n \geq 0} a_n u^n) = \sum_{n \geq 0} a_n\otimes u^n, ~~p_2(\sum_{n \geq 0} a_n u^n) = \sum_{n \geq 0} 1\otimes a_n u^n, ~~a_n \in \mathcal{O}_{L_0},
\]
which extends to $p_1, p_2\colon S \rightarrow S(1)$. Since $\mathcal{O}_{L_0}(1)$ is $p$-torsion free, the maps $p_1, p_2\colon \mathcal{O}_{L_0} \rightarrow \mathcal{O}_{L_0}(1)$ are flat. So $p_1, p_2\colon S/(p^n) \rightarrow S(1)/(p^n)$ are flat for each $n \geq 1$ by \cite[Tag~07HD]{stacks-project}.

By crystalline base change, we have isomorphisms
\[
H^i_{\mathrm{cris}}(\mathcal{X}_{k'} / \mathcal{O}_{L_0})\otimes_{\mathcal{O}_{L_0}, p_1} \mathcal{O}_{L_0}(1) \cong H^i_{\mathrm{cris}}(\mathcal{X}_{k'} / \mathcal{O}_{L_0}(1)) \cong H^i_{\mathrm{cris}}(\mathcal{X}_{k'} / \mathcal{O}_{L_0})\otimes_{\mathcal{O}_{L_0}, p_2} \mathcal{O}_{L_0}(1)
\]
and
\[
H^i_{\mathrm{cris}}(\mathcal{X}_0 / S)\otimes_{S, p_1} S(1) \cong H^i_{\mathrm{cris}}(\mathcal{X}_0 / S(1)) \cong H^i_{\mathrm{cris}}(\mathcal{X}_0 / S)\otimes_{S, p_2} S(1).
\]
By the standard constuction, these naturally give connections
\[
\nabla_{M^i}\colon M^i \rightarrow M^i\otimes_{\mathcal{O}_{L_0}} \widehat{\Omega}_{\mathcal{O}_{L_0}}, ~~\nabla_{\mathscr{M}^i}\colon \mathscr{M}^i \rightarrow \mathscr{M}^i\otimes_{\mathcal{O}_{L_0}} \widehat{\Omega}_{\mathcal{O}_{L_0}}.
\]
Furthermore, it follows from above constructions that $\nabla_{\mathscr{M}^i}$ is compatible with $\nabla_{\mathscr{D}(M^i[p^{-1}])}$ via the isomorphism $\mathscr{M^i}[p^{-1}] \cong \mathscr{D}(M^i[p^{-1}]) \coloneqq S\otimes_{\mathcal{O}_{L_0}} M^i[p^{-1}]$ given in Proposition \ref{prop:comparison-crys-cohom's}. 

Let $\widetilde{S}(1)$ be the $p$-adically completed PD-envelope of $\mathcal{O}_{L_0}[x]\widehat{\otimes}_{W(k)}\mathcal{O}_{L_0}[y]$ with respect to the kernel of the map $\mathcal{O}_{L_0}[x]\widehat{\otimes}_{W(k)}\mathcal{O}_{L_0}[y] \rightarrow \mathcal{O}_L/(p)$ given by $x, y \mapsto \pi$. Note that the kernel is generated by $p$, $E(x)$, $x-y$, and $a\otimes 1 - 1\otimes a, ~a \in \mathcal{O}_{L_0}$. Similarly as above, we have the maps $p_1, p_2\colon S \rightarrow \widetilde{S}(1)$ given by $p_1(u) = x$ and $p_2(u) = y$. For each $n \geq 1$, the induced map $p_1, p_2\colon S/(p^n) \rightarrow \widetilde{S}(1)/(p^n)$ is flat by \cite[Tag~07HD]{stacks-project}, since the $p$-adically completed PD-envelope of $\mathcal{O}_{L_0}(1)[z]$ with respect to $(z)$ (for $z = x-y$) is flat over $\mathcal{O}_{L_0}$. So by crystalline base change, we have isomorphisms 
\begin{equation} \label{eq:crys-cohom-HPD}
H^i_{\mathrm{cris}}(\mathcal{X}_0 / S)\otimes_{S, p_1} \widetilde{S}(1) \cong H^i_{\mathrm{cris}}(\mathcal{X}_0 / \widetilde{S}(1)) \cong H^i_{\mathrm{cris}}(\mathcal{X}_0 / S)\otimes_{S, p_2} \widetilde{S}(1).	
\end{equation}

For $\mathfrak{S} = \mathcal{O}_{L_0}[\![u]\!]$, write $\widehat{\Omega}_{\mathfrak{S}}$ for the $p$-adic completion of $\Omega^1_{\mathfrak{S} / \mathbf{Z}_p}$. By \cite[Prop.~1.3.1]{berthelot-messing-dieudonne-cryst-III}, we have $\widehat{\Omega}_{\mathfrak{S}} \cong (\bigoplus_{j=1}^d \mathfrak{S}\cdot dT_j) \oplus \mathfrak{S}\cdot du$. The universal derivation $d\colon \mathfrak{S} \rightarrow \widehat{\Omega}_{\mathfrak{S}}$ naturally extends to $d\colon S \rightarrow S\otimes_{\mathfrak{S}} \widehat{\Omega}_{\mathfrak{S}}$. By \cite[Thm.~6.6]{berthelot-ogus-book}, the HPD-stratification given by the isomorphisms in (\ref{eq:crys-cohom-HPD}) induces a connection
\[
\nabla\colon \mathscr{M}^i \rightarrow \mathscr{M}^i\otimes_{\mathfrak{S}} \widehat{\Omega}_{\mathfrak{S}} \cong (\bigoplus_{j=1}^d \mathscr{M}^i \cdot dT_j)\oplus \mathscr{M}^i \cdot du,
\]
which is compatible with $\nabla_{\mathscr{M}^i}$ above. For $1 \leq j \leq d$, let $\partial_{T_j}\colon \mathscr{M}^i \rightarrow \mathscr{M}^i$ be the derivation given by $\nabla$ composed with the projection to the $dT_j$-component, and let $\partial_u\colon \mathscr{M}^i \rightarrow \mathscr{M}^i$ be the derivation given by $\nabla$ composed with the projection to the $du$-component. Write $N_{T_j} \coloneqq T_j\partial_{T_j}$ and $N \coloneqq -u\partial_u$ for the corresponding derivations on $\mathscr{M}^i$. 

\begin{prop} \label{prop:rational-cryst-comparison}
We have a natural isomorphism
\[
H^i_{\text{\'et}}(X_{\overline{L}}, \mathbf{Q}_p)\otimes_{\mathbf{Q}_p} \mathbf{OB}_{\mathrm{cris}}(\mathcal{O}_{\overline{L}}) \cong H^i_{\mathrm{cris}}(\mathcal{X}_{k'} / \mathcal{O}_{L_0})\otimes_{\mathcal{O}_{L_0}}\mathbf{OB}_{\mathrm{cris}}(\mathcal{O}_{\overline{L}})
\]	
compatible with filtration, $\varphi, \nabla$, and $G_L$-actions.
\end{prop}

\begin{proof}
Let $\mathbf{C}_p\coloneqq \widehat{\overline{L}} = \widehat{\overline{K_g}}$, and let $\mathcal{O}_{\mathbf{C}_p}$ its ring of integers. By \cite[Thm.~1.1]{bhatt-morrow-scholze-integralpadic}, since $H^i_{\mathrm{cris}}(\mathcal{X}_{k'} / \mathcal{O}_{L_0})\otimes_{\mathcal{O}_{L_0}} W(k_g) \cong H^i_{\mathrm{cris}}(\mathcal{X}_{k_g} / W(k_g))$, we have a natural isomorphism
\[
\alpha\colon H^i_{\text{\'et}}(X_{\mathbf{C}_p}, \mathbf{Q}_p)\otimes_{\mathbf{Q}_p} \mathbf{B}_{\mathrm{cris}}(\mathcal{O}_{\mathbf{C}_p}) \cong H^i_{\mathrm{cris}}(\mathcal{X}_{k'} / \mathcal{O}_{L_0})\otimes_{\mathcal{O}_{L_0}}\mathbf{B}_{\mathrm{cris}}(\mathcal{O}_{\mathbf{C}_p})
\] 	
compatible with filtration, $\varphi$, and $G_{K_g}$-actions. Consider the associated isomorphism
\[
\beta\colon H^i_{\text{\'et}}(X_{\mathbf{C}_p}, \mathbf{Q}_p)\otimes_{\mathbf{Q}_p} \mathbf{B}_{\mathrm{cris}}(\mathcal{O}_{\mathbf{C}_p}) \cong H^i_{\mathrm{cris}}(\mathcal{X}_{\mathcal{O}_{\mathbf{C}_p}/p} / \mathbf{A}_{\mathrm{cris}}(\mathcal{O}_{\mathbf{C}_p}))\otimes_{\mathbf{A}_{\mathrm{cris}}(\mathcal{O}_{\mathbf{C}_p})} \mathbf{B}_{\mathrm{cris}}(\mathcal{O}_{\mathbf{C}_p}) 
\] 
as in \cite[Thm.~1.8]{bhatt-morrow-scholze-integralpadic}. Since $\beta$ is functorial on $\mathcal{X}$ which is defined over $\mathcal{O}_L$, it is compatible with $G_L$-actions. Note that
\[
H^i_{\mathrm{cris}}(\mathcal{X}_{\mathcal{O}_{\mathbf{C}_p}/p} / \mathbf{A}_{\mathrm{cris}}(\mathcal{O}_{\mathbf{C}_p})) \cong H^i_{\mathrm{cris}}(\mathcal{X}_0 / S)\otimes_S \mathbf{A}_{\mathrm{cris}}(\mathcal{O}_{\mathbf{C}_p})
\]
by crystalline base change, since $S/(p^n) \rightarrow \mathbf{A}_{\mathrm{cris}}(\mathcal{O}_{\mathbf{C}_p}))/(p^n)$ is flat for each $n \geq 1$. Let $\sigma \in G_L$, and consider two maps $S \rightarrow \mathbf{A}_{\mathrm{cris}}(\mathcal{O}_{\mathbf{C}_p})$ and the composite $S \rightarrow \mathbf{A}_{\mathrm{cris}}(\mathcal{O}_{\mathbf{C}_p}) \stackrel{\sigma}{\rightarrow}\mathbf{A}_{\mathrm{cris}}(\mathcal{O}_{\mathbf{C}_p})$. Via $p_1, p_2\colon S \rightarrow \widetilde{S}(1)$, these two maps induce $h_{\sigma}\colon \widetilde{S}(1) \rightarrow \mathbf{A}_{\mathrm{cris}}(\mathcal{O}_{\mathbf{C}_p})$. Then $\sigma$-action on $H^i_{\mathrm{cris}}(\mathcal{X}_{\mathcal{O}_{\mathbf{C}_p}/p} / \mathbf{A}_{\mathrm{cris}}(\mathcal{O}_{\mathbf{C}_p}))$ is induced by 
\[
H^i_{\mathrm{cris}}(\mathcal{X}_0 / \widetilde{S}(1))\otimes_{\widetilde{S}(1), h_{\sigma}} \mathbf{A}_{\mathrm{cris}}(\mathcal{O}_{\mathbf{C}_p}) \cong H^i_{\mathrm{cris}}(\mathcal{X}_{\mathcal{O}_{\mathbf{C}_p}/p} / \mathbf{A}_{\mathrm{cris}}(\mathcal{O}_{\mathbf{C}_p}))
\] 
together with the isomorphisms
\[
H^i_{\mathrm{cris}}(\mathcal{X}_0 / S)\otimes_{S, p_1} \widetilde{S}(1) \cong H^i_{\mathrm{cris}}(\mathcal{X}_0 / \widetilde{S}(1)) \cong H^i_{\mathrm{cris}}(\mathcal{X}_0 / S)\otimes_{S, p_2} \widetilde{S}(1).
\]
Thus, $\sigma$-action on $\mathscr{M}^i\otimes_S \mathbf{A}_{\mathrm{cris}}(\mathcal{O}_{\mathbf{C}_p})$ is given by Equation (\ref{eq:Galois-action}), with $N$ and $N_{T_j}$'s induced from the isomorphisms (\ref{eq:crys-cohom-HPD}) as above.

By Proposition \ref{prop:comparison-crys-cohom's}, above implies that $G_L$-action on $M^i[p^{-1}]\otimes_{L_0}\mathbf{B}_{\mathrm{cris}}(\mathcal{O}_{\mathbf{C}_p}) = M^i[p^{-1}]\otimes_{L_0, \iota_2}\mathbf{B}_{\mathrm{cris}}(\mathcal{O}_{\overline{L}})$ is given by Equation (\ref{eq:galois-action-D}). On the other hand, consider the projection 
\[
pr\colon M^i[p^{-1}]\otimes_{L_0, \iota_1} \mathbf{OB}_{\mathrm{cris}}(\mathcal{O}_{\overline{L}}) \rightarrow M^i[p^{-1}] \otimes_{L_0, \iota_2} \mathbf{B}_{\mathrm{cris}}(\mathcal{O}_{\overline{L}})
\]
as in Section \ref{sec:strongly-div-latt}. Equation (\ref{eq:section}) induces 
\[
s\colon M^i[p^{-1}] \otimes_{L_0, \iota_2} \mathbf{B}_{\mathrm{cris}}(\mathcal{O}_{\overline{L}}) \cong (M^i[p^{-1}]\otimes_{L_0, \iota_1} \mathbf{OB}_{\mathrm{cris}}(\mathcal{O}_{\overline{L}}))^{\nabla = 0}
\]   
compatible with filtration, $\varphi$, and $G_L$-actions, where the $G_L$-action on $M^i[p^{-1}]\otimes_{L_0, \iota_1} \mathbf{OB}_{\mathrm{cris}}(\mathcal{O}_{\overline{L}})$ is given by the trivial action on $M^i[p^{-1}]$. Thus, $\alpha$ gives an isomorphism
\[
\alpha\colon T^i[p^{-1}]\otimes_{\mathbf{Q}_p}\mathbf{B}_{\mathrm{cris}}(\mathcal{O}_{\overline{L}}) \cong (M^i[p^{-1}]\otimes_{L_0, \iota_1} \mathbf{OB}_{\mathrm{cris}}(\mathcal{O}_{\overline{L}}))^{\nabla = 0}
\]
compatible with filtration, $\varphi$, and $G_L$-actions.

Note that $M^i[p^{-1}] \in \mathrm{MF}(\varphi, \nabla)$ is weakly admissible since $M_i[p^{-1}]\otimes_{L_0} W(k_g)[p^{-1}] = H^i_{\mathrm{cris}}(\mathcal{X}_{k_g} / W(k_g))[p^{-1}]$ is weakly admissible. From the isomorphism $\alpha$, we obtain $V_{\mathrm{cris}}(M^{i}[p^{-1}]) = (T^i[p^{-1}])^{\vee}$. Hence, we deduce from Theorem \ref{thm:weakly-adm-cris} that
\[
H^i_{\text{\'et}}(X_{\overline{L}}, \mathbf{Q}_p)\otimes_{\mathbf{Q}_p} \mathbf{OB}_{\mathrm{cris}}(\mathcal{O}_{\overline{L}}) \cong H^i_{\mathrm{cris}}(\mathcal{X}_{k'} / \mathcal{O}_{L_0})\otimes_{\mathcal{O}_{L_0}, \iota_1}\mathbf{OB}_{\mathrm{cris}}(\mathcal{O}_{\overline{L}})
\]	
compatibly with filtration, $\varphi, \nabla$, and $G_L$-actions.
\end{proof}

By above proposition and its proof, $T^i[p^{-1}]$ is a crystalline representation of $G_L$ with Hodge-Tate weights in $[-i, 0]$, and we have $D_{\mathrm{cris}}(T^i[p^{-1}]^\vee) \cong M^i[p^{-1}]$. By Proposition \ref{prop:comparison-crys-cohom's}, we can identify $\mathscr{M}^i[p^{-1}] = \mathscr{D}(M^i[p^{-1}]) = S[p^{-1}]\otimes_{L_0} M^i[p^{-1}]$. 

\begin{prop} \label{prop:crys-cohom-strly-div-latt}
Under the conditions of \emph{Theorem \ref{thm:cryst-cohom-strly-div-latt}}, $\mathscr{M}^i$ is a strongly divisible lattice in $\mathscr{D}(D^i)$.
\end{prop}

\begin{proof}
Note that $H^j_{\mathrm{cris}}(\mathcal{X}_{k_g} / W(k_g))$ for $j = i, ~i+1$ are $p$-torsion free. So by \cite[Thm.~5.4]{cais-liu-breuil-kisin-mod-cryst-cohom}, $\mathscr{M}^i_g \coloneqq H^i_{\mathrm{cris}}(\mathcal{X}_{\mathcal{O}_{K_g}, 0} / S_g)$ is finite free over $S_g$ of rank $m \coloneqq \mathrm{rk}_{W(k_g)[p^{-1}]} D_{\mathrm{cris}, K_g}(T^i[p^{-1}]^\vee) = \mathrm{rk}_{L_0} D_{\mathrm{cris}}(T^i[p^{-1}]^\vee)$.

Since $b_g\colon \mathcal{O}_{L_0}[\![u]\!] \rightarrow W(k_g)[\![u]\!]$ is flat, the induced map $S/(p^n) \rightarrow S_g/(p^n)$ is flat for each integer $n \geq 1$ by \cite[Tag~07HD]{stacks-project}. Moreover, $S/(p^n) \rightarrow S_g/(p^n)$ is a map of local rings, so it is faithfully flat. By crystalline base change, we have
\[
\mathscr{M}^i/p^n\mathscr{M}^i \otimes_S S_g \cong \mathscr{M}^i_g/p^n\mathscr{M}^i_g.
\]
In particular, $\mathscr{M}^i/p\mathscr{M}^i$ is free over $S/(p)$ of rank $m$ by faithfully flat descent. Let $\{e_1, \ldots, e_m\}$ be a basis for $\mathscr{M}^i/p\mathscr{M}^i$ over $S/(p)$, and choose a lifting $\hat{e}_1, \ldots, \hat{e}_m \in \mathscr{M}^i$. By Nakayama's lemma, $\mathscr{M}^i$ is generated by $\hat{e}_1, \ldots, \hat{e}_d$ as $S$-module. Since $\mathscr{M}^i[p^{-1}]$ is free over $S[p^{-1}]$ of rank $m$ by Proposition \ref{prop:comparison-crys-cohom's}, we conclude that $\mathscr{M}^i$ is free over $S$ of rank $m$. Furthermore, the map $\mathscr{M}^i\otimes_S S_g \rightarrow \mathscr{M}^i_g$ is an isomorphism.

Let $\mathrm{Fil}^i \mathscr{M}^i = \mathscr{M}^i \cap \mathrm{Fil}^i \mathscr{D}(M^i[p^{-1}])$. Since $\mathscr{M}^i_g$ is a strongly divisible lattice of weight $i$ in $S_g\otimes_{W(k_g)} D_{\mathrm{cris}, K_g}(T^i[p^{-1}]^\vee)$ by \cite[Thm.~5.4]{cais-liu-breuil-kisin-mod-cryst-cohom}, we have
\[
\varphi(\mathrm{Fil}^i \mathscr{M}^i_g) \subset p^i \mathscr{M}^i_g.
\]
Since $\mathrm{Fil}^i \mathscr{M} \subset \mathrm{Fil}^i \mathscr{M}_g$ and $S \cap p^i S_g = p^i S$,
\[
\varphi(\mathrm{Fil}^i \mathscr{M}) \subset \mathscr{M}^i \cap p^i \mathscr{M}^i_g = p^i\mathscr{M}^i.
\]

It remains to study the differential operators on $\mathscr{M}^i$. We observed above that $\mathscr{M}^i$ is stable under $\nabla_{\mathscr{D}(M^i[p^{-1}])}$. Furthermore, since $\mathscr{M}^i_g$ is a strongly divisible lattice in $S_g\otimes_{W(k_g)} D_{\mathrm{cris}, K_g}(T^i[p^{-1}]^\vee)$, we have
\[
N_{\mathscr{D}(M^i[p^{-1}])}(\mathscr{M}^i) \subset \mathscr{M}_g^i \cap \mathscr{D}(M^i[p^{-1}]) = \mathscr{M}^i.
\]
\end{proof}

\begin{proof}[Proof of Theorem \ref{thm:cryst-cohom-strly-div-latt}]
It remains to show that 	$T_{\mathrm{cris}}(\mathscr{M}^i) \cong (T^{i})^{\vee}$ as $G_L$-representations. Let $\mathscr{N}^i \in \mathrm{Mod}_S^{i}$ such that $T_{\mathrm{cris}}(\mathscr{N}^i) = (T^{i})^{\vee}$ given by Theorem \ref{thm:equiv-category-strly-div-latt-cris-reps}. Note that $\mathscr{N}^i$ is constructed explicitly in Section \ref{sec:essential-surjectivity}, compatibly with the construction in \cite{liu-semistable-lattice-breuil} via $b_g\colon S \rightarrow S_g$. We identify $\mathscr{M}^i[p^{-1}] = \mathscr{D}(M^i[p^{-1}]) = \mathscr{N}^i[p^{-1}]$.

Since $\mathscr{N}^i_g \coloneqq \mathscr{N}^i\otimes_S S_g$ is the strongly divisible lattice of $S_g\otimes_{W(k_g)} D_{\mathrm{cris}, K_g}(T^i[p^{-1}]^\vee)$ such that $T_{\mathrm{cris}}(\mathscr{N}^i_g) = (T^{i})^{\vee}$ as $G_{K_g}$-representations, it is shown in \cite[Pf. of Thm. 5.4]{cais-liu-breuil-kisin-mod-cryst-cohom} that $\mathscr{M}^i_g = \mathscr{N}^i_g$ as $S_g$-submodules of $S_g\otimes_{W(k_g)} D_{\mathrm{cris}, K_g}(T^i[p^{-1}]^\vee) = \mathscr{D}(M^i[p^{-1}])\otimes_S S_g$. Thus,
\[
\mathscr{M}^i = \mathscr{M}^i[p^{-1}] \cap \mathscr{M}^i_g = \mathscr{N}^i[p^{-1}] \cap \mathscr{N}^i_g = \mathscr{N}^i.
\]
\end{proof}

\bibliographystyle{amsplain}
\bibliography{library}
	
\end{document}